\documentclass[11pt,a4paper]{article}
\usepackage{amsfonts,amsmath,amsthm,epsf,graphics,verbatim,amssymb,amscd,eucal,bbm}    

\newtheorem{theorem}{Theorem}[section]
\newtheorem{corollary}[theorem]{Corollary}
\newtheorem{proposition}[theorem]{Proposition}
\newtheorem{lemma}[theorem]{Lemma}
\theoremstyle{definition}

\newtheorem{definition}[theorem]{Definition}

\theoremstyle{remark}

\def\Blackboardfont{\mathbb}

\newcommand{\moins}{ {\setminus} }
\newcommand{\set}[2]{\{#1\mid\;#2\}}
\newcommand{\pres}[2]{\langle \: #1 \mid #2 \: \rangle}

\def\Z{{\Blackboardfont Z}}

\def\N{{\Blackboardfont N}}
\def\R{{\Blackboardfont R}}

\def\cB{{\mathcal B}}

\def\cS{{\mathcal S}}

\def\cX{{\mathcal X}}

\def\iff{\Longleftrightarrow}
\def\eref#1{(\ref{#1})}

\begin{document}
\sloppy

\title{\bf Random walks on free products of cyclic groups}

\author{Jean {\sc Mairesse}
\thanks{LIAFA, CNRS-Universit\'e Paris 7, case
    7014, 2, place Jussieu, 75251 Paris Cedex 05, France. E-mail: {\tt
      Jean.Mairesse@liafa.jussieu.fr}} \ and Fr\'ed\'eric {\sc Math\'eus} 
\thanks{LMAM, Universit\'e de Bretagne-Sud, Campus de Tohannic,  BP
  573, 56017 Vannes, France. E-mail: {\tt Frederic.Matheus@univ-ubs.fr}}}

\maketitle

\begin{abstract}
Let $G$ be a free product of a finite family of finite groups,
with the set of generators being formed by the union of the finite
groups. We consider a transient nearest-neighbor random walk on $G$. 
We give a new proof of the fact that the  harmonic measure is
a special Markovian measure entirely determined 
by a finite set of polynomial equations. 
We show that in 
several simple cases of interest, the polynomial equations
can be explicitly solved, to get closed form formulas for the
drift. The examples considered are  
$\Z/2\Z\star \Z/3\Z$, $\Z/3\Z\star \Z/3\Z$,  
$\Z/k\Z\star \Z/k\Z$, and the Hecke groups $\Z/2\Z\star \Z/k\Z$.
We also use these various examples to study Vershik's notion of extremal
generators, which is based on the relation between the drift,
the entropy, and the volume of the group. 
\end{abstract}

\smallskip


\smallskip

\textsl{AMS classification (2000):} Primary 60J10, 60B15, 60J22, 65C40;
Secondary 28D20; 37M25.


\section{Introduction}
\label{se-intro} 

The properties of the harmonic measure associated with
a nearest-neighbor random walk on a free group, or a free product of groups, have been studied
by many authors, see \cite{DyMa,SaSt,lall93,ledr00,woes86b}, or the
monograph \cite{woes} and
the references therein. 
In this context, the Green kernel has a multiplicative
structure. Consequently, the harmonic measure is Markovian, see in
particular \cite[Section 5]{SaSt} and \cite[Section 6]{woes86b}. 
In \cite{mair04}, this property is viewed from a different angle. 
It is proved that the harmonic measure is a Markovian measure with a
special combinatorial structure, called {\em Markovian
  multiplicative}. It is entirely determined by 
its initial distribution, 
which is itself characterized as the unique solution to a finite set
of polynomial equations coined as the {\em Traffic Equations}. 
The result of \cite{mair04} is proved for a whole class 
of pairs formed by a group (or monoid) and a finite set of generators:
the so-called {\em 0-automatic pairs}.
The property of being Markovian multiplicative is very specific.
For instance, in the related context of trees with finitely many cone types,
the harmonic measure is Markovian but not
Markovian multiplicative, see \cite[Section 5]{NaWo}. 


\medskip

In this paper, we focus on nearest-neighbor random walks (NNRW) on free products of finite
groups, and in particular of finite cyclic groups.
There may be several motivations for specifically studying such random
walks. First, they are among the simplest non-commutative
random walks. As such, they serve as a reference point, and numerous
results have first been proved in this context before being extended.
Second, they may be pertinent in the 
physics of polymers. This is discussed and argued in \cite{nech,NeVo}.

\medskip

There are three types of results being proved in the paper. 

\medskip

{\em Section \ref{se-fpmhm}} - We revisit the result of \cite{mair04}.
Consider a free product of finite groups $G$, the set of
generators $\Sigma$ being the union of the finite groups (the {\em natural}
generators).  The pair $(G,\Sigma)$ is a special case of 
0-automatic pair. 
Consider a transient NNRW on
$(G,\Sigma)$, say $(X_n)_n$. 
Assume that the
group elements are written in reduced form. 
The {\em harmonic measure} is the law of $X_{\infty}=\lim_n X_n$.
First, we give a short proof of the special Markovian
structure of the harmonic measure. This proof is different from the
one in \cite{mair04} and would not work in the more general context of
0-automatic pairs. Second, we take advantage of the restricted setting
to prove more precise results. In particular, we characterize the
cases where the harmonic measure is stationary with respect to the 
translation shift. 

\medskip

{\em Section \ref{se-expli}} - The result on the harmonic measure has
interesting computational 
consequences. Indeed, in many situations, 
the Traffic Equations can be solved 
``explicitly'', in order to get a closed form formula
for the drift, the entropy or the minimal harmonic functions.
We illustrate this by explicitly computing the drift in the following 
cases: (i) the general NNRW 
on the modular group $\Z/2\Z\star \Z/3\Z$, (ii) a two parametrized families
of NNRW on $\Z/3\Z\star \Z/3\Z$, and (iii) the simple (with respect
to minimal generators) NNRW on 
$\Z/k\Z\star \Z/k\Z$, and on the Hecke groups $\Z/2\Z\star \Z/k\Z$. 

\medskip

{\em Section \ref{se-vershik}} - We investigate Vershik notion of extremal
generators, which is 
based on the link between drift, entropy, and volume for random walks
on a group~\cite{vers00}. We prove the following. For
a free product of finite groups, the set of natural generators
is extremal. This uses the
special structure of the harmonic measure (see Section \ref{se-fpmhm}). We also show that
for the group $\Z/4\Z\star \Z/4\Z$, the minimal set of generators is not
extremal. This is achieved by explicitly computing the drift and the
entropy (as in Section \ref{se-expli}). 

\medskip

In our view, the interest of the present work is to provide 
a large collection of tractable models,
whereas very few were previously available. This is a potential source
for examples or counter-examples.
In particular, none of the drift computations in the paper appeared in the
litterature before. For the 
few examples of non-elementary explicit computations previously
available, see \cite{derr,DyMa,ledr00,NaWo,NeVo,SaSt}. 
Nevertheless, there exists an alternative potential method for the
effective computation of the drift (not the entropy)
which is due to Sawyer and Steger \cite{SaSt}.
In this approach, the drift is expressed as a functional of the {\em first-passage
generating series} of the random walk. The simplest of our
computational results can also be retrieved using this method.
We detail and discuss this approach in \S
\ref{sse-other}.
The Sawyer-Steger method links
the problem of computing the drift with the problem of
computing the generating series of transition probabilities. Concerning the
latter problem, there exists an important litterature, especially for
random walks on free groups and free products, see \cite{cart88,CaSo86,GeWo,woes86} and
\cite[Sections II.9 and III.17]{woes}.

\medskip

For much more material on random walks on discrete infinite groups (including
aspects not even touched upon here, e.g., 
boundary theory, or central/local limit theorems), see
\cite{woes,kaim00,vers00} and the references there. 

\medskip

Part of the results presented in this paper, as well as in the companion
paper \cite{mair04}, were announced without proofs in the Proceedings of
the International Colloquium of Mathematics and Computer Science held
in Vienna in September 2004~\cite{MaMa04}.

\section{Preliminaries}\label{se-prel}

\paragraph{Notations.}

Let $\N$ be the set of non-negative integers and $\N^*$ the set of
positive integers.
We denote the support of a random variable by $\text{supp}$. 
If $\mu$
is a measure on a group $(G,\ast)$, then 
$\mu^{*n}$ is the $n$-fold convolution product of $\mu$, that is the
image of the product measure $\mu^{\otimes n}$ by the product map
$G\times \cdots \times G \rightarrow G, \ (g_1,\ldots, g_k)\mapsto
g_1\ast g_2\ast \cdots \ast g_k$.
The symbol $\sqcup$ is used for the disjoint union of sets.  
Given a finite set $\Sigma$, a
vector $x\in \R^{\Sigma}$, and $S\subset \Sigma$, set 
$x(S)=\sum_{u\in S} x(u)$. 

\subsection{Random walks on groups}

Consider a finitely generated group $(G,\ast)$ with unit element $1_G$.
Let $\Sigma\subset G$ be a finite set of generators of $G$ (with $1_G\not\in \Sigma$ and
$u\in \Sigma \implies u^{-1}\in \Sigma$). 
The {\em length} with respect to $\Sigma$
of a group element $u$ is:
\begin{equation}\label{eq-length}
|u|_{\Sigma}= \min \set{k}{u=s_1\ast \cdots \ast s_k, s_i \in \Sigma}\:.
\end{equation}

The {\em Cayley graph} $\cal X$$(G,\Sigma)$ of a group $G$ with respect to a set of 
generators $\Sigma$ is the graph with $G$ as set of vertices and with an edge 
between $u$ and $v$ if and only if $u^{-1}v \in \Sigma$. 
Observe that $|u|_{\Sigma}$ 
is the geodesic distance from $1_G$ to $u$ in the Cayley graph. 
 
\medskip

Let $\mu$ be a probability distribution
over $\Sigma$. Consider the Markov chain on the state space $G$ with one-step transition
probabilities given by: $\forall g\in G, \forall a\in \Sigma, \
P_{g,g\ast a}=\mu(a)$. This Markov chain is called the {\em random walk}
  (associated with) $(G,\mu)$. 
  It is a {\em nearest neighbor} random walk: one-step moves occur 
  between nearest neighbors in the Cayley
  graph $\cX(G,\Sigma)$. When $\mu(s)=1/|\Sigma|$ for all $s\in \Sigma$, we say that
  the random walk is {\em simple}. 

Let $(x_n)_n$ be a sequence of
i.i.d. r.v.'s distributed according to $\mu$. Set
\begin{equation}\label{eq-rw}
X_0=1, \ X_{n+1}=X_n \ast x_n=x_0\ast x_1\ast \cdots \ast x_n \:.
\end{equation}
The sequence $(X_n)_n$ is a {\em realization} of the random walk
$(G,\mu)$. The law of $X_n$ is $\mu^{*n}$.
Since $|u\ast v|_{\Sigma} \leq |u|_{\Sigma} +|v|_{\Sigma}$, Guivarc'h \cite{guiv80} observed that a simple
corollary of Kingman's Subadditive Ergodic Theorem \cite{king} 
 is the existence of a constant $\gamma \in \R_+$ such that a.s.
 and in $L^p$, for all $1\leq p <\infty$,
\begin{equation}\label{eq-drift}
\lim_{n\rightarrow \infty} \frac{|X_n|_{\Sigma}}{n} = \gamma \:.
\end{equation}
We call $\gamma$ the {\em drift}. Intuitively, $\gamma$ is the speed of escape to infinity
of the walk. 

\subsection{Free products and harmonic measure}

Let $(G_i)_{i\in I}$ be a finite family of finite groups, with
$|I|\geq 2$. Let $1_{G_i}$ be
the unit of $G_i$. Set $\Sigma_i=G_i\moins \{1_{G_i}\}$ and 
set $\Sigma=\sqcup_i \Sigma_i$.
Let $\iota:\Sigma \rightarrow I$ be defined by 
$\iota(u)= j$ if $u\in \Sigma_j$. It is also convenient to set
$\Sigma_a=\Sigma_{\iota(a)}$ for all $a\in \Sigma$. 

Let $\Sigma^*$ be the free monoid over the alphabet $\Sigma$ and denote its
unit, the empty word, by 1. 
Define the set of {\em normal form words} 
$L\subset \Sigma^*$ by 
\begin{equation}\label{eq-nfw}
L= \{u_1\cdots u_k \in \Sigma^*, \ \forall i \in \{1,\cdots, k-1\},
\iota(u_i)\neq \iota(u_{i+1})\}\:. 
\end{equation}
Hence, $L$ consists of all words over the alphabet $\Sigma$ whose consecutive letters come from
different subalphabets $\Sigma_i$. 
Observe that $1\in L$. 

The {\em free product} $G=\star_{i\in I} G_i$ is
the group with set of elements $L$, unit element 1, and group law
$\ast$ defined recursively by: 
\[
u_1\cdots u_k \ast v_1 \cdots v_l = \begin{cases}
u_1\cdots (u_{k-1})(u_k)(v_1)(v_2)\cdots v_l & \mbox{if } \iota(u_k)\neq
\iota(v_1) \\
u_1\cdots (u_{k-1})(u_k \ast v_1)(v_2)\cdots v_l & \mbox{if } \iota(u_k) = 
\iota(v_1), \ u_k \neq v_1^{-1} \\
u_1\cdots u_{k-1} \ast v_2 \cdots v_l & \mbox{if }  u_k = v_1^{-1} \end{cases}\:,
\]
where in the second case, $(u_k\ast v_1)$ is the product in $G_{\iota(u_k)}$
of $u_k$ and $v_1$. Roughly, the law of $G$ is the concatenation with
possible simplifications at the contact point to reach a normal form word. 

The {\em length} of an element $u$ of $(\star_{i\in I} G_i)$ is the
length (i.e. number of letters) of the word $u$ in $L$. We denote it
by $|u|$. Observe that $|u|=\min \{k \mid 
u_1\ast \cdots \ast u_k = u, \ u_i\in \Sigma \}= |u|_{\Sigma}$. 

\medskip

Let $\mu$ be a probability measure over $\Sigma$ such that
$\bigcup_{n\in \N^*} \mbox{supp} \ \mu^{*n} = \star_{i\in I} G_i$.
Let $(X_n)_n$ be a realization of the random walk $(\star_{i\in I}
G_i,\mu)$ as defined above. The sequence $(X_n)_n$ can be viewed as a Markov chain
on $L$.
Below the drift is defined according to \eref{eq-drift}
with respect to the length $|\cdot|$. 

\medskip

The group $\Z/2\Z\star \Z/2\Z$ is amenable, and 
any nearest neighbor random walk on it is recurrent. 
But apart from this group,
all the free products considered are
non-amenable. 
Therefore, if $G$ is not $\Z/2\Z\star \Z/2\Z$, any
random walk living on the whole group is 
transient and has a strictly positive drift (see \cite{guiv80} and 
\cite[Chapter 1.B]{woes} for details). 

>From now on, the random walks considered are assumed to be
transient. Equivalently, we work on a free product group $G$ different
from $\Z/2\Z\star \Z/2\Z$ and the support of $\mu$ generates the whole group. 

\medskip

Consider the set $\Sigma^{\N}$ equipped with the product topology. 
Denote by $(u_1\cdots u_n \Sigma^{\N})$ the order-$n$ cylinder in $\Sigma^\N$
defined by $u_1\cdots u_n$. 
Define the set of {\em (right) infinite normal form words}
$L^{\infty}\subset \Sigma^{\N}$ by 
\begin{equation}\label{eq-Linfty}
L^{\infty}= \{u_0u_1u_2\cdots   \in \Sigma^{\N}, \ \forall i \in \N^*,
\iota(u_i)\neq \iota(u_{i+1})\}\:. 
\end{equation}
A word belongs to $L^{\infty}$ iff all its finite prefixes belong to
$L$. 

Consider the map $\Sigma\times L^{\infty} \rightarrow
L^{\infty}, (a,\xi) \mapsto a\cdot \xi$, with $a\cdot
\xi=a\xi_0\xi_1\cdots$ if $\iota(a)\neq \iota(\xi_0)$, $a\cdot
\xi=(a\ast \xi_0)\xi_1\cdots$ if $\iota(a)= \iota(\xi_0), a\neq
\xi_0^{-1}$, and  
 $a\cdot \xi=\xi_1\xi_2\cdots$ if $a=\xi_0^{-1}$. 
Equip $\Sigma^{\N}$ with the Borel $\sigma$-algebra associated with
the product topology. This induces a $\sigma$-algebra on
$L^{\infty}$. 
Given a measure $\nu^{\infty}$ on $L^{\infty}$ and $a\in \Sigma$,
define the measure $a\nu^{\infty}$ by: $\int f(\xi) d(a\nu^{\infty})(\xi) = \int
f(a\cdot \xi) d\nu^{\infty}(\xi)$.
A probability measure $\nu^{\infty}$ on $L^{\infty}$ is {\em
$\mu$-invariant} if 
\begin{equation}\label{eq-inva}
\nu^{\infty}(\cdot)  =  
\sum_{a\in \Sigma} \mu(a) [a\nu^{\infty}] (\cdot) \:.
\end{equation}

\begin{proposition}\label{pr-harm}
There exists a r.v. $X_{\infty}$ valued in $L^{\infty}$ such that a.s.
\[
\lim_{n\rightarrow \infty} X_n = X_{\infty}\:,
\]
meaning that the length of the common prefix between
$X_n$ and $X_{\infty}$ goes to infinity a.s. Let $\mu^{\infty}$ be the
distribution of $X^{\infty}$. The probability $\mu^{\infty}$ is $\mu$-invariant
and is the only $\mu$-invariant probability  on $L^{\infty}$. We call it the {\em harmonic
  measure} of $(G,\mu)$. The drift of the random walk
is given by:
\begin{eqnarray}\label{eq-drift2}
\gamma  & = & \sum_{a\in \Sigma} \mu(a) \Bigl[ -\mu^{\infty}(a^{-1}\Sigma^{\N}) +
  \sum_{b\in \Sigma\moins \Sigma_a} \mu^{\infty}(b\Sigma^{\N})\Bigr] \:.
\end{eqnarray}
In words, $\gamma$ is the expected change of length of an infinite
normal form distributed according to $\mu^{\infty}$, when
left-multiplied by an element distributed according to $\mu$. 
\end{proposition}

In the context of the free group, this is proved for instance in
\cite[Theorem 1.12, Theorem 4.10]{ledr00}. The proofs adapt easily to the present
setting. Several of the key arguments go back to \cite{furs,furs71}, see
\cite{ledr00} for precise references. 

Intuitively, the harmonic measure $\mu^{\infty}$ gives the direction in which $(X_n)_n$
goes to infinity. 

\section{Free Products Have a Markovian Harmonic
  Measure}\label{se-fpmhm}

In \cite[Theorems 4.5 and 5.3]{mair04}, it is proved that the harmonic
measure for random walks on 0-automatic pairs has a special Markovian
multiplicative 
structure. In this section, we revisit the result. 
We concentrate on a subclass of 0-automatic pairs: the pairs formed by
free products of finite groups with natural generators. In this
setting, we get a more elementary proof of the result. We also refine
the result by discussing the cases where the harmonic measure is
shift-invariant (on top of being $\mu$-invariant). To that purpose, we
associate two sets of equations with the random walk: the {\em
  Traffic Equations} (as in \cite{mair04}) but also the {\em Stationary
  Traffic Equations}. The results in Propositions \ref{pr-stat} and
\ref{pr-two} are new. 

\bigskip

Define 
$\mathring{\cB}= \{ x \in \R^{\Sigma} \mid \forall u, \: x(u) > 0, \ \sum_u x(u) =1 \}$.
Consider $r\in \mathring{\cB}$. 
Define the matrix $P$ of dimension $\Sigma\times \Sigma$ by 
\begin{equation}\label{eq-transi}
P_{u,v}=\begin{cases} r(v)/r(\Sigma\moins \Sigma_u) &
                      \text{if } v\in \Sigma\moins \Sigma_u \\ 
                      0    & \text{otherwise } 
\end{cases}\:.
\end{equation}
It is the transition matrix of an irreducible Markov
Chain on the state space $\Sigma$. Set $p=(r(a)r(\Sigma\moins
\Sigma_a), a\in \Sigma)$ and 
$\pi=p/p(\Sigma)$. Observe that $\pi P=\pi$. In words, 
$\pi$ is the stationary distribution of the Markov chain defined by
$P$. 

Let $(U_n)_n$ be a realization of the Markov chain with transition
matrix $P$ and starting from $U_1$ such that $P\{U_1=x\}=r(x)$. Set
$U^{\infty}=\lim_n U_1\cdots U_n$, and let $\nu^{\infty}$ be the distribution
of $U^{\infty}$. 
Clearly the support of $\nu^{\infty}$ is included in $L^{\infty}$. For 
$u_1\cdots u_k \in L$, we have
\begin{eqnarray}\label{eq-mult}
\nu^{\infty}(u_1\cdots u_k \Sigma^{\N}) & = & r(u_1)P_{u_1,u_2}\cdots P_{u_{k-1},u_k}
\\ 
&=& r(u_1)
\frac{r(u_2)}{r(\Sigma\moins \Sigma_{u_1})} \cdots
\frac{r(u_k)}{r(\Sigma\moins \Sigma_{u_{k-1}})}
\ = \ \frac{r(u_1)}{r(\Sigma\moins \Sigma_{u_1})}
   \cdots 
\frac{r(u_{k-1})}{r(\Sigma\moins \Sigma_{u_{k-1}})} r(u_k)\:.\nonumber
\end{eqnarray} 
We call $\nu^{\infty}$ the {\em Markovian multiplicative} probability
measure associated with $r$.   

\medskip

The measure $\nu^{\infty}$ is in
general non-stationary with respect to the translation shift $\tau:
\Sigma^{\N}\rightarrow \Sigma^{\N}, (x_n)_n \mapsto
(x_{n+1})_n$. Indeed, the distribution of the first marginal 
is $r$ which is different in general from the stationary distribution
$\pi$. 
An important special case is when $\nu^{\infty}$ is nevertheless
stationary and ergodic, i.e. when $r=\pi$. This happens if and only if 
\begin{equation}\label{eq-stationary}
\forall i \in I, \quad r(\Sigma_i)= 1/|I|\:.
\end{equation}

\begin{definition}[Traffic Equations]
The {\em Traffic Equations} associated with $(G,\mu)$ are defined by:
$\forall a\in \Sigma$, 
\begin{equation}\label{eq-traffic}
x(a) = \mu(a)\sum_{u\in \Sigma \moins \Sigma_a} x(u) + \sum_{u\ast v = a}
\mu(u)x(v) + \sum_{u\in \Sigma \moins \Sigma_a} \mu(u^{-1})
\frac{x(u)}{\sum_{v\in \Sigma\moins \Sigma_u} x(v)} x(a) \:.
\end{equation}
\end{definition}

The Traffic Equations are closely related to the harmonic measure of
$(G,\mu)$. Next lemma is proved in a more general context in
\cite[Lemma 5.2]{mair04}.

\begin{lemma}\label{le-traffic} 
If the harmonic measure $\mu^{\infty}$ is the Markovian multiplicative
measure associated with $r\in \mathring{\cB}$, then $r$ is a solution to the Traffic
Equations \eref{eq-traffic}. Conversely, if the Traffic Equations admit a
solution $r\in  \mathring{\cB}$, then
the harmonic measure $\mu^{\infty}$ is the Markovian multiplicative 
measure associated with $r$.  
\end{lemma}

Using \eref{eq-stationary}, we can complete the statement of Lemma
\ref{le-traffic} as follows. 

\begin{definition}[Stationary Traffic Equations]
The {\em Stationary Traffic Equations} associated with $(G,\mu)$ are defined by:
$\forall a\in \Sigma$, 
\begin{equation}
\label{eq-trafficstat}
x(a) = \mu(a)\frac{|I|-1}{|I|} +  \sum_{u\ast v =a} \mu(u) x(v) + 
x(a)\frac{|I|}{|I|-1} \sum_{u\in \Sigma\moins \Sigma_a} \mu(u^{-1}) x(u) \:.
\end{equation}
\end{definition}

\begin{lemma}\label{le-trafficstat} 
The harmonic measure $\mu^{\infty}$ is Markovian multiplicative
associated with $r$ 
{\em and} ergodic if and only if the Stationary Traffic Equations
admit a solution $r$ in $\mathring{\cB}$. 
\end{lemma}

A corollary of Lemma \ref{le-traffic}, resp. Lemma
\ref{le-trafficstat}, 
is that the Traffic Equations, resp. the
Stationary Traffic Equations, have at most one solution in
$\mathring{\cB}$. The Stationary Traffic Equations do not always have
 solution. But the Traffic Equations do, see Theorem
 \ref{th-jackpot}. 

\begin{theorem}\label{th-jackpot}
Let $G=\star_{i\in I} G_i$ be the free product of a finite family of finite groups, with
$|I|\geq 2$, and $\forall i, |G_i|>1$.   Let
$\mu$ be a probability measure on $\Sigma=\sqcup_i G_i\moins \{1_{G_i}\}$. 
Assume that $\bigcup_{n\in \N^*} \text{supp} \ \mu^{*n} =
G$ and that the random walk
$(G,\mu)$ is transient. Then the Traffic Equations \eref{eq-traffic}
have a unique solution $r\in \mathring{\cB}$.
The harmonic measure of the random walk
is the Markovian multiplicative measure associated with $r$. 
\end{theorem}

For finitely generated free groups, 
the special Markovian structure of the harmonic measure
is a classical result \cite{DyMa,SaSt,ledr00}.
Theorem \ref{th-jackpot} is proved in \cite{mair04} in the more
general context of 0-automatic pairs.
Here we give a short proof of Theorem \ref{th-jackpot} which is close
in spirit to the proofs in \cite{DyMa,SaSt,ledr00,NaWo}. On the other
hand, the proof below is quite different from the one in \cite{mair04}
and would not work in the general context of 0-automatic pairs. 

\begin{proof}
For all $a\in \Sigma$, define $q(a)=P \{ \exists n \mid X_n
=a\}$, the probability of ever
hitting $a$. Clearly, $0<q(a)<1$. Besides, we have:
$\forall a\in \Sigma$, 
\begin{equation}\label{eq-reach}
q(a) = \mu(a) + \sum_{u \ast v = a} \mu(u) q(v) + 
q(a)\sum_{c\in \Sigma\moins \Sigma_a} \mu(c)q(c^{-1})\:. 
\end{equation}
The first two terms on the right-hand side of \eref{eq-reach} are more
or less obvious. Now assume that the random walk starts with an
initial step of type $c\in \Sigma\moins\Sigma_a$. Given the tree-like
structure of the Cayley graph, it has to go back to 1 before possibly
reaching $a$. Now, the probability of ever hitting 1 starting from $c$
is equal to the probability of ever hitting $c^{-1}$ starting from
1. This accounts for the third right-hand term in \eref{eq-reach}. 

A simple rewriting of the Traffic Equations \eref{eq-traffic} gives: 
\begin{equation}\label{eq-traffic2}
\frac{x(a)}{x(\Sigma\moins\Sigma_a)}  =  \mu(a) + \sum_{u\ast v =a}
\mu(u) \frac{x(v)}{x(\Sigma\moins\Sigma_v)}  + \frac{x(a)}{x(\Sigma\moins\Sigma_a)} 
\sum_{u\in \Sigma\moins \Sigma_a} \mu(u^{-1}) 
\frac{x(u)}{x(\Sigma\moins\Sigma_u)} \:.
\end{equation}
Hence it is natural to look for a solution $r$ to the Traffic Equations
satisfying 
\begin{equation}\label{eq-natural}
\forall a\in \Sigma, \quad
\frac{r(a)}{r(\Sigma\moins\Sigma_a)}=q(a)\:.
\end{equation}
It remains to be proved that the Equations \eref{eq-natural} have a solution in
$r$. Clearly, they have a solution iff the following equations have a
solution:
\begin{equation}\label{eq-eq}
\begin{cases}
r(\Sigma_i)=q(\Sigma_i)r(\Sigma\moins \Sigma_i)  & (i) \\
\sum_{j\in I} r(\Sigma_j) =1   & (\text{sum})
\end{cases}\:.
\end{equation}
For a given $i$, consider
the Equations: $\{(i),(\text{sum})\}$. The solution is 
$r(\Sigma_i) = q(\Sigma_i)/(1+q(\Sigma_i)), \ r(\Sigma\moins \Sigma_i) =
1/(1+q(\Sigma_i))$. 
In order to have a global solution to \eref{eq-eq}, the necessary and
sufficient condition is that:
\begin{equation}\label{eq-consistent}
\sum_{i\in I}   \frac{q(\Sigma_i)}{1+q(\Sigma_i)} =1 \:.
\end{equation}
Now let us prove \eref{eq-consistent}. 
Starting from \eref{eq-reach} and summing over $\Sigma_i$, we get:
\begin{eqnarray*}
q(\Sigma_i) &  = & \mu(\Sigma_i) + \sum_{a\in \Sigma_i} \sum_{u\ast v=a}
\mu(u)q(v) + q(\Sigma_i)\sum_{a\in \Sigma\moins \Sigma_{i}}
\mu(a^{-1}) q(a) \\
 &  = &  \mu(\Sigma_i)+ \mu(\Sigma_i)q(\Sigma_i) -\sum_{a\in \Sigma_i}
\mu(a^{-1})q(a)  + q(\Sigma_i)\sum_{a\in \Sigma\moins \Sigma_{i}}
\mu(a^{-1}) q(a) 
\end{eqnarray*}
It follows that 
\begin{eqnarray*}
\frac{q(\Sigma_i)}{1+q(\Sigma_i)} &  = &   \mu(\Sigma_i) +
\frac{q(\Sigma_i)}{1+q(\Sigma_i)} \sum_{a\in \Sigma\moins \Sigma_{i}}
\mu(a^{-1}) q(a) - \frac{1}{1+q(\Sigma_i)} \sum_{a\in \Sigma_{i}}
\mu(a^{-1}) q(a) \\
\sum_i \frac{q(\Sigma_i)}{1+q(\Sigma_i)}  &  = & 1  + \sum_i \bigl[\sum_{a\in
    \Sigma_i} \mu(a^{-1})q(a) \bigr] \bigl[\sum_{j\neq i} \frac{q(\Sigma_j)}{1+q(\Sigma_j)} -
    \frac{1}{1+q(\Sigma_i)} \bigr] \\
\sum_i \frac{q(\Sigma_i)}{1+q(\Sigma_i)} - 1 &  = & \bigl[ \sum_{a\in
    \Sigma} \mu(a^{-1})q(a) \bigr] \bigl[ \sum_i
    \frac{q(\Sigma_i)}{1+q(\Sigma_i)} - 1 \bigr] \:.
\end{eqnarray*}
Since $0<q(a)<1$ for all $a$, we have  $\sum_{a\in
    \Sigma} \mu(a^{-1})q(a)<1$. We conclude that we must have: $ \sum_i
    q(\Sigma_i)/(1+q(\Sigma_i)) - 1 =0$. Hence Equation
    \eref{eq-consistent} holds. It implies that the Traffic Equations
    have a solution which is:
\begin{equation}\label{eq-inverse}
\forall a\in   \Sigma, \quad r(a) = \frac{q(a)}{1+q(\Sigma_a)} \:.
\end{equation}
According to Lemma \ref{le-traffic}, such a solution is necessarily unique, and the
harmonic measure is the Markovian multiplicative measure associated
with $r$. 
\end{proof}

\begin{corollary}\label{co-drift}
Under the assumptions of Theorem \ref{th-jackpot}, the drift is given by
\begin{equation}\label{eq-drift3}
\gamma  =  \sum_{a\in \Sigma} \mu(a)\bigl[ -r(a^{-1}) + \sum_{b \in
    \Sigma\moins \Sigma_a}
r(b) \bigr] \:,
\end{equation}
where $r$ is the unique solution in $\mathring{\cB}$ to the Traffic
Equations. 
\end{corollary}

\begin{corollary}\label{co-harmo}
Under the assumptions of Theorem \ref{th-jackpot}, we have 
\[
\forall a \in \Sigma, \qquad P\{\exists n \mid X_n =a\} =
r(a)/r(\Sigma\moins \Sigma_a)\:,
\]
where $r$ is the unique solution in $\mathring{\cB}$ to the Traffic
Equations. 
\end{corollary}


It follows from Lemma \ref{le-trafficstat} and Theorem
\ref{th-jackpot} that the harmonic measure is shift-invariant iff the
Stationary Traffic Equations have a solution. In Proposition
\ref{pr-stat}, we give a sufficient condition for this to happen.

\begin{proposition}\label{pr-stat}
Let $H$ be a finite group and let $(G_i)_{i\in I}$ be a finite family
of copies of $H$. Let $\pi_i$ be the
isomorphism between $G_i$ and $H$. Let $\nu$ be a probability measure
on $H\moins\{1_H\}$. 
Consider the free product $G=\star_{i\in I} G_i$
and let $\mu$ be the
probability measure on $\Sigma=\sqcup_i G_i\moins \{1_{G_i}\}$ defined by:
$\forall g \in G_i\moins \{1_{G_i}\}, \ \mu(g)=\nu\circ \pi_i(g)/|I|$.  
Then the harmonic measure of $(G,\mu)$ is stationary and ergodic with
respect to the translation shift $\tau: 
\Sigma^{\N}\rightarrow \Sigma^{\N}, (x_n)_n \mapsto
(x_{n+1})_n$.
\end{proposition}

\begin{proof}
Set $\Sigma_i = G_i\moins \{1_{G_i}\}$ for all $i$. 
Let $a\in \Sigma_i$ and $b\in \Sigma_j$ be such that
$\pi_i(a)=\pi_j(b)$. By a symmetry argument, we have
$r(a)=r(b)$, where $r$ is the solution to the Traffic Equations. 
A direct consequence is that $r(\Sigma_i)=r(\Sigma_j)$ for all
$i,j$. Hence Condition \eref{eq-stationary} is satisfied and 
the harmonic measure is ergodic. 
\end{proof}

\begin{proposition}\label{pr-two}
Consider a random walk $(G,\mu)$ where $G=G_1\star G_2$ is the
free product of two arbitrary finite groups. 
Then the harmonic measure $\mu^{\infty}$ satisfies: 
$\forall u \in L, \forall k \in \N, \
\mu^{\infty}(u\Sigma^{\N})=
\mu^{\infty}(\Sigma^{2k}u\Sigma^{\N})$. That is, $\mu^{\infty}$ 
is stationary and ergodic with
respect to the shift $\tau^2: 
\Sigma^{\N}\rightarrow \Sigma^{\N}, (x_n)_n \mapsto
(x_{n+2})_n$.
\end{proposition}

\begin{proof}
Applying \eref{eq-consistent}, we get:
\begin{equation}\label{eq-q1q2=1}
\frac{q(\Sigma_1)}{1+q(\Sigma_1)} + \frac{q(\Sigma_2)}{1+q(\Sigma_2)}
=1 \ \implies \ q(\Sigma_1)q(\Sigma_2)=1 \:.
\end{equation}
Consider $u=u_1\cdots \in L$ with for instance $u_1\in \Sigma_1$. We
have, using \eref{eq-traffic}:
\begin{eqnarray*}
\mu^{\infty}(\Sigma^{2}u\Sigma^{\N}) & = & \sum_{v_1\in \Sigma_1}\sum_{v_2\in
  \Sigma_2} \mu^{\infty}(v_1v_2u\Sigma^{\N}) 
\ = \ \sum_{v_1\in \Sigma_1} q(v_1)\sum_{v_2\in
  \Sigma_2} q(v_2)\mu^{\infty}(u\Sigma^{\N})  \\
& = & q(\Sigma_1)q(\Sigma_2) \mu^{\infty}(u\Sigma^{\N}) \ = \
  \mu^{\infty}(u\Sigma^{\N})\:.
\end{eqnarray*}
We prove in the same way that: $\forall u \in L, \forall k \in \N, \
\mu^{\infty}(u\Sigma^{\N})= \mu^{\infty}(\Sigma^{2k}u\Sigma^{\N})$.
\end{proof}

There exists no simple analog of Proposition \ref{pr-two} for the free
product of three or more finite groups. 

\medskip

The Identity \eref{eq-q1q2=1}: $q(\Sigma_1)q(\Sigma_2)=1$, is quite
unexpected. 
Indeed it can be rephrased as: the average number of different elements
visited in $\Sigma_1$ is the
inverse of the average number of different elements visited in $\Sigma_2$.

\paragraph{Free products of countable groups} $ $ \medskip

Consider a free product $G=\star_{i\in I} G_i$ of a finite family of {\em countable} groups.
Let $\mu$ be a probability measure on the countable set $\Sigma=\sqcup_i G_i\moins \{1_{G_i}\}$. 
Assume that $\bigcup_{n\in \N^*} \text{supp} \ \mu^{*n} = G$ and that the random walk
$(G,\mu)$ is transient. One can define the set of Traffic Equations
(resp. Stationary Traffic Equations) 
exactly as in \eref{eq-traffic} (resp. \eref{eq-trafficstat}). It is a set of infinitely
many equations involving infinite sums. 
We have the following. 

\begin{theorem}
The statements of Lemmas \ref{le-traffic} and \ref{le-trafficstat},
 Theorem \ref{th-jackpot}, Corollaries \ref{co-drift} 
 and \ref{co-harmo}, Propositions \ref{pr-stat} and \ref{pr-two}, remain true for a free product 
$G=\star_{i\in I} G_i$ of a finite family of countable groups.
\end{theorem} 

\begin{proof}
The only difference with the case of finite groups is that we have to prove that
$q(\Sigma_i)$ is finite for all $i\in I$. We proceed as follows.
For $u\in G$, define the r.v. $\tau(u)=\min\{n \mid X_n =u \}$ (with
$\tau(u)=\infty$ if $u$ is not reached). Define the series $q(u,z) = \sum_{n\in \N}
P\{\tau(u)=n\} z^n$. Observe that $q(u,1)=q(u)$, the probability of
ever hitting $u$. The family of series $(q(u,z))_{u\in \Sigma}$ satisfies
the following version of \eref{eq-reach}:

\begin{equation}\label{eq-reach-ser}
\forall a\in \Sigma,\:
q(a,z) = z\,\mu(a) + z\sum_{u \ast v = a} \mu(u) q(v,z) + 
z\,q(a,z)\sum_{c\in \Sigma\moins \Sigma_a} \mu(c)q(c^{-1},z)\:.
\end{equation}

For $i\in I$, we set
$q(\Sigma_i,z)=\sum_{u\in \Sigma_i} q(u,z)$, $R_{i}(z)=\sum_{u\in \Sigma_i}
\mu(u^{-1})q(u,z)$, and  $S_{i}(z)=\sum_{j\neq i}R_{j}(z)$. 
These series are well-defined since all the coefficients are bounded
by 1. Observe that, since $q(a)<1$, we have $R_{i}(1)<\mu(\Sigma_i)$
and therefore $S_{i}(1)<1-\mu(\Sigma_i)$.
\par
Now, starting from \eref{eq-reach-ser} and summing over $\Sigma_i$, we get:
\[
q(\Sigma_i,z)
\:=\:z\,\mu(\Sigma_i)+z\,\mu(\Sigma_i)q(\Sigma_i,z)-z\,R_{i}(z)+z\,q(\Sigma_i,z)S_{i}(z) \:.
\]
Therefore $q(\Sigma_i)=q(\Sigma_i,1)=(\mu(\Sigma_i)-R_{i}(1))/(1-\mu(\Sigma_i)-S_{i}(1))$ is finite.
\end{proof}

As an example of this situation, consider the
free group $G=\Z\star\Z$. Denote by $a$ and $b$ the generators of the
two factors, so that $\Sigma=\{a^{i},i\in \Z\setminus\{0\}\}\cup\{b^{j},j\in \Z\setminus\{0\}\}$.
Let $\mu$ be a probability measure on $\Sigma$ whose support generates
$G$. The harmonic measure $\mu^{\infty}$ of the random walk
$(G,\mu)$ satisfies, e.g., $\mu^{\infty}(a^{k_1}b^{k_2}\cdots
a^{k_{\ell}}\Sigma^{\N}) = q(a^{k_1})q(b^{k_2})\cdots
q(b^{k_{\ell}-1}) r(a^{k_{\ell}})$. When $\mu$ is
concentrated on $\{a,a^{-1},b,b^{-1}\}$, the following
simplifications hold (see \cite{mair04}): $\forall k>0, \
q(a^k)=q(a)^k, \ r(a^k)=q(a)^{k-1}r(a), 
\ \forall k<0, \ q(a^k)=q(a^{-1})^{-k}, \ r(a^k)=
q(a^{-1})^{-k-1}r(a^{-1})$, and the analog for $b$. 

\section{Explicit Drift Computations}
\label{se-expli}

In Theorem \ref{th-jackpot}, the harmonic measure is completely
determined via the vector $r$ which is itself the solution of an
explicit finite set of polynomial equations of degree 2. 
In small or simple examples, it is possible to go further, that is, to
solve these equations to get closed form formulas for the harmonic
measure, and therefore the drift. It is the
program that we now carry out. We compute the drift for
several specific and interesting 
cases of free products of two cyclic groups. 

\medskip

The details of the computations and complete proofs of the results are
not given. They can be found in an appendix 
posted on the Math ArXiv~\cite{MaMa05}. 
Some of the results have been obtained with the help of Maple and
Mathematica. 

\subsection{Comparison with other methods for computing the
  drift}\label{sse-other}

We first discuss alternative existing methods for computing the drift.
(They do not work for computing the entropy for instance.)

\medskip

Let $G=G_1\star G_2$ be a free product of two finite groups. Set
$\Sigma_i= G_i\moins\{1_{G_i}\}$ and
$\Sigma = \Sigma_1 \sqcup \Sigma_2$. Let $\mu$ be a
probability measure on $\Sigma$ such that:
$\forall i, \forall x \in \Sigma_i, \ \mu(x)=\mu(\Sigma_i)/\#
\Sigma_i$.
In words, $\mu$ is uniform on each of the two groups.
Consider the random walk $(G,\mu)$. Here, computing
the drift becomes elementary and does not require knowing that the
harmonic measure is Markovian.

Set $p=\mu(\Sigma_1), k_1=\#\Sigma_1,$ and $k_2=\# \Sigma_2$.
Denote by $i\in \{1,2\}$, the set of elements
of $G$ whose normal form representative ends with a letter in $\Sigma_i$.
When we are far from the unit element $1_G$, the random walk $(X_n)_n$
on $G$ induces a Markov chain on $\{1,2\}$ with
transition matrix:
\[
P=\left[ \begin{array}{cc}  p(k_1 -1)/k_1 & p/k_1 + 1-p \\ (1-p)/k_2
    +p & (1-p)(k_2 -1)/k_2 \end{array}\right]\:.  
\]
Let $\pi$ be the stationary distribution, that is $\pi P=\pi,
\pi(1)+\pi(2)=1$.
By
the Ergodic Theorem for Markov Chains, we have $\lim_n [\ P\{X_n\in
  1\} \ ,\ P\{X_n\in 2\} \ ]=\pi$. The value of the drift follows readily:
\begin{equation}\label{eq-elementary-general}
\gamma  =  \lim_n \frac{1}{n} \sum_{i=0}^{n-1}
E[|X_{i+1}|_{\Sigma}- |X_i|_{\Sigma}] + \frac{|X_0|_{\Sigma}}{n} =
E_{\pi}[|X_1|_{\Sigma}-|X_0|_{\Sigma}]
 =  \frac{2p(1-p)(k_1k_2-1)}{(1-p)k_1+pk_2 +k_1k_2} \:. 
\end{equation}
Now assume that $G=G_1\star \cdots \star G_k$, where the $G_i$ are finite
groups, and assume that $\forall i, \forall x \in \Sigma_i=G_i\moins
\{1_{G_i}\}, \ \mu(x)=\mu(\Sigma_i)/\# 
\Sigma_i$. Then each of the finite groups can be collapsed into a
single node, and the random walk $(G,\mu)$ projects into a nearest
neighbor randow walk on a tree with $k$ cone types, using the
terminology of \cite{NaWo}. 
In particular the
formulas for the drift given in \cite{NaWo} apply. 

\medskip

None of the formulas obtained in \S \ref{sse-z2z3}-\S \ref{sse-z2zk}
correspond to the above two situations.

\medskip

Now let us discuss the Sawyer-Steger method \cite{SaSt}. It was developped for the
free group but adapts to the present situation.
Let $G=G_1\star \cdots \star G_k$ be a free product of finite groups. Set
$\Sigma_i= G_i\moins\{1_{G_i}\}$ and $\Sigma = \sqcup_i \Sigma_i$.
For $g\in G$, define the r.v. $\tau(g)=\min\{n \mid X_n =g \}$ (with
$\tau(g)=\infty$ if $g$ is not reached). Define the {\em first-passage
  generating series} $S\in \R[[y,z]]$ by:
\begin{equation}\label{eq-sast1}
S(y,z) = \sum_{k\in \N} y^k \sum_{|g|_{\Sigma} =k} \ \sum_{n\in \N}
P\{\tau(g)=n\} z^n \:.
\end{equation}
Let $S_y$ and $S_z$ denote the partial derivatives of $S$ with respect
to $y$ and $z$.
Adapting the results in \cite[Theorem 2.2 and Section 6]{SaSt} (see
also \cite[Section 6]{NaWo}), one obtains the following formula for the
drift:
\begin{equation}\label{eq-sast2}
\gamma = S_y(1,1)/S_z(1,1)\:.
\end{equation}
For $u\in G$, define the series $q(u,z) = \sum_{n\in \N}
P\{\tau(u)=n\} z^n$. Observe that $q(u,1)=q(u)$, the probability of
ever hitting $u$, defined at the beginning of the proof of Theorem \ref{th-jackpot}.
In particular, if one encapsulates the Equations \eref{eq-reach} as
$q(a)=\Phi_a (q)$, then $q(a,z)=z\Phi_a(q(z))$. Using this last set of
Equations, as well as the corresponding set of Equations for the derivatives $dq(a,z)/dz$,
and playing around with the Equations \eref{eq-sast1} and \eref{eq-sast2}, we get:
\begin{equation}\label{eq-sast3}
\gamma = \frac{ \sum_{i=1}^k q_i/(1+q_i)^2 } {
  \sum_{i=1}^k q_i' /(1+q_i)^2 } \:, \qquad q_i=\sum_{u\in \Sigma_i}
  q(u,1), \quad q_i'=
\sum_{u\in \Sigma_i} \bigl[ \frac{dq(u,z)}{dz}\bigr]_{|z=1}\:.
\end{equation}

This formula is more complicated than the one obtained in
\eref{eq-drift3}. Our approach centers around the knowledge
that $\mu^{\infty}$ is Markovian multiplicative, which is a more
direct path. Consequently,
it gives more chances to solve the equations
to get a closed form formula. As an exercise, we tried to retrieve the
results for $\gamma$ in \S \ref{sse-z2z3}-\S \ref{sse-z2zk} using
\eref{eq-sast3}. We succeeded in two cases: formulas \eref{eq-driftz2z3}
and \eref{eq-driftz3z3a=b}. On the other hand, the
results in \S \ref{sse-simplezn}-\S \ref{sse-z2zk} seem totally out of
reach.

\subsection{Random walks on $\Z/2\Z\star \Z/3\Z$}\label{sse-z2z3}

The group
$\Z/2\Z\star \Z/3\Z$  is isomorphic to the modular group PSL$(2,\Z)$,
i.e. the group of 2x2 matrices with integer entries and determinant
$1$, quotiented by $\pm$Id. 

\begin{figure}[ht]
\[ \epsfxsize=220pt \epsfbox{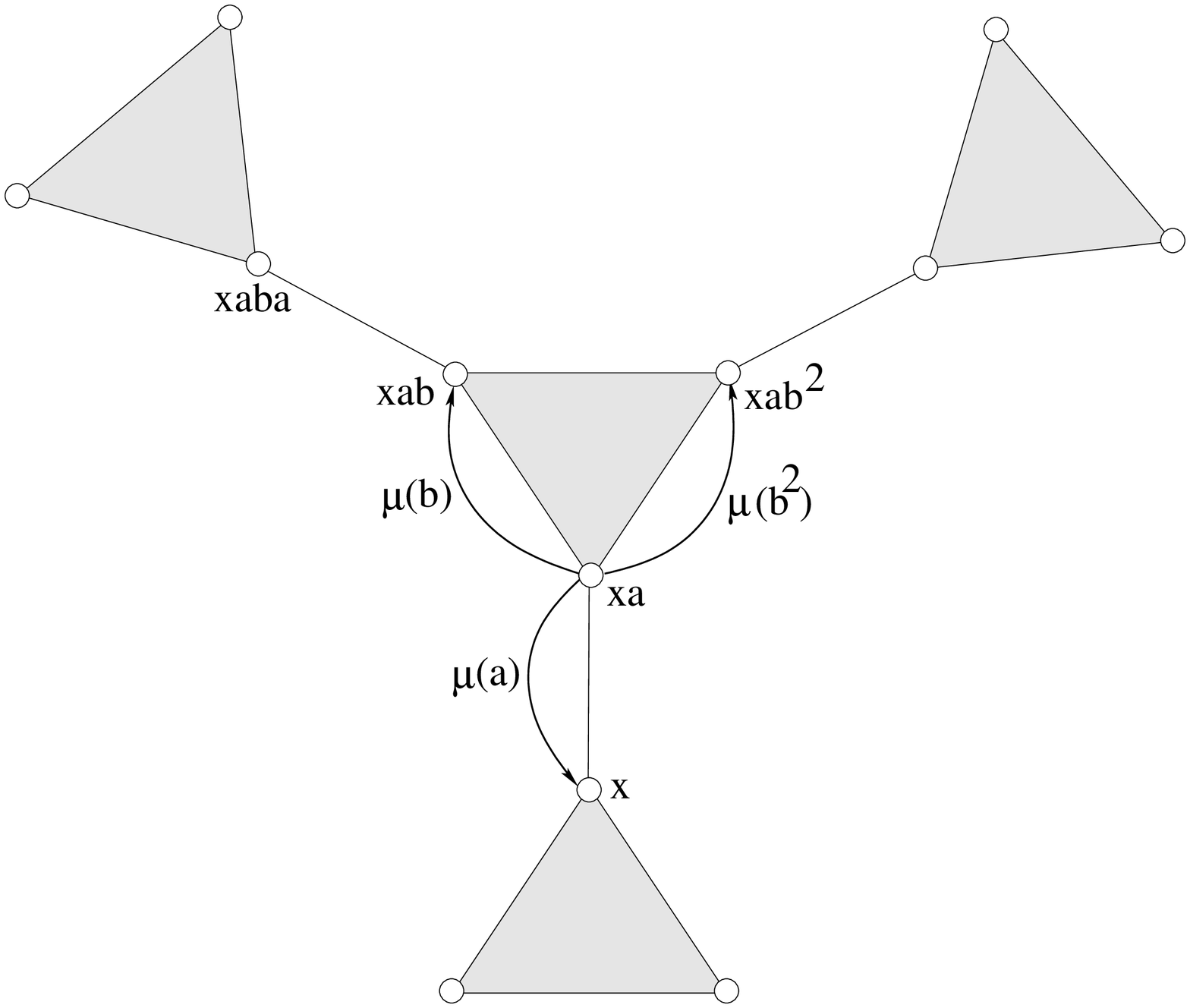} \]
\caption{A nearest neighbor random walk on $\Z/2\Z\star \Z/3\Z$.}
\label{fi-z2z3} 
\end{figure}

Consider a general nearest neighbor random walk $(\Z/2\Z\star
\Z/3\Z,\mu)$. Set $\mu(a)=1-p-q, \mu(b)=p,\mu(b^2)=q$. In
Figure \ref{fi-z2z3}, we have represented the Cayley graph of
$\Z/2\Z\star \Z/3\Z$ and the one-step transitions of the random walk. 

\begin{figure}[ht]
\[ \epsfxsize=190pt \epsfbox{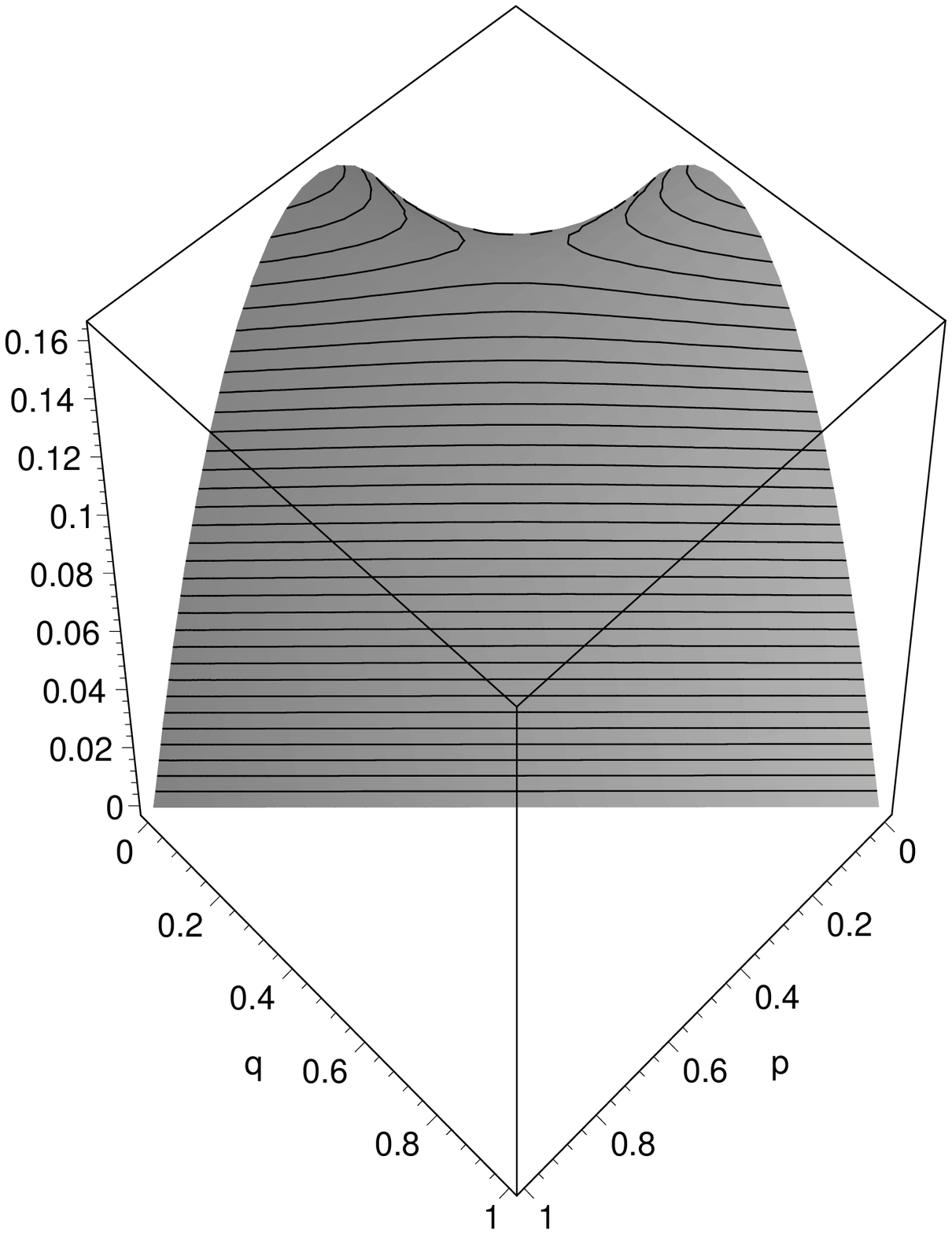} \qquad
 \epsfxsize=190pt \epsfbox{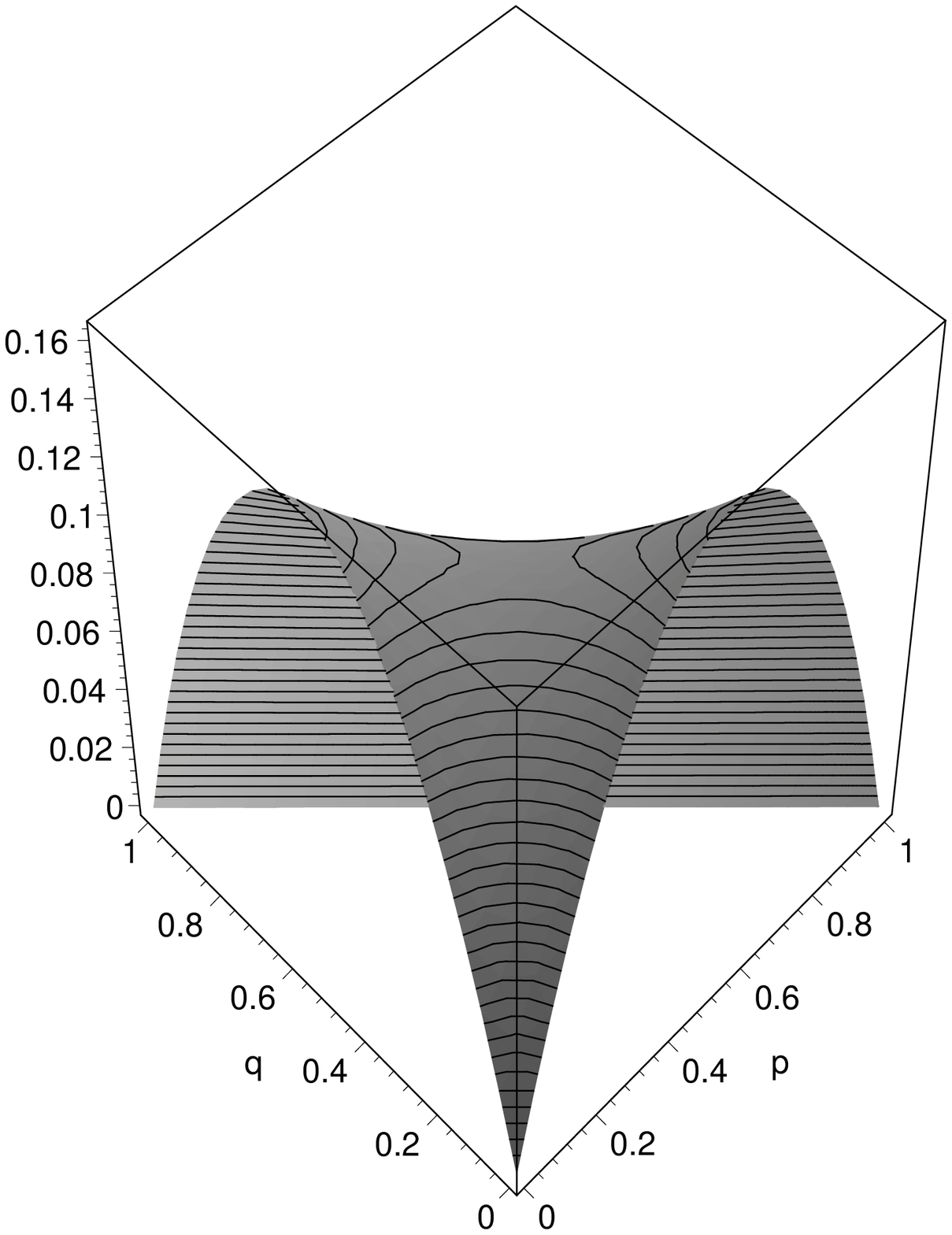} \]
\caption{The drift of the random walk $(\Z/2\Z\star \Z/3\Z,\mu)$ as a
  function of $p=\mu(b)$ and $q=\mu(b^2)$, from two different angles.}
\label{fi-z2z3drift} 
\end{figure}

Using Theorem \ref{th-jackpot} and solving the Traffic Equations
explicitly, we get~:

\begin{eqnarray*}
r(a) & = & \frac{p^2+q^2-2pq-p-q+4 - \Delta_1}{2\Delta_2} \\
r(b) & = & \frac{q^3-3q^2 +p^2q -5pq+2p+6q -
  (2-q)\Delta_1}{2(q-p)\Delta_2} \\
r(b^2) & = & \frac{p^3 -3p^2+pq^2 -5pq +6p + 2q -
  (2-p)\Delta_1}{2(p-q)\Delta_2}  \:,
\end{eqnarray*}
with
\begin{eqnarray*}
\Delta_1 & = &
\sqrt{p^4+q^4-2p^3-2q^3+2p^2q^2-6p^2q-6pq^2+5p^2+5q^2+6pq} \\
\Delta_2 & = & p^2+q^2-pq-2p-2q+4 \:.
\end{eqnarray*}
Set $r=\mu(a)=1-p-q$. By Corollary \ref{co-drift}, the drift is
then: 
\begin{equation}\label{eq-driftz2z3}
\gamma =  \frac{2r\Bigl( pq-p-q+
  \sqrt{(p^2+q^2)(3+(r+p)^2+(r+q)^2)+2pq(2r+1)}  
\Bigr)}{(r+p)^2+(r+q)^2 -pq+2  } \:.
\end{equation}
As a curiosity, the drift is maximized for $r=z_0,p=1-z_0,q=0$ (or
$r=z_0,p=0,q=1-z_0)$), where $z_0$ is the root of
$[z^6+12z^4-4z^3+47z^2-48z+12]$ whose numerical value is
$0.490275\cdots$. The corresponding numerical value of the drift is
$\gamma_{\max}= 0.163379\cdots$. This was not a priori obvious!

\medskip

We are in the domain of application of Proposition \ref{pr-two}. On
the other hand,  
the harmonic measure is never shift-invariant.
This can be proved by showing directly that the Stationary
Traffic Equations have no solution. 


\subsection{Random walks on $\Z/3\Z\star \Z/3\Z$}\label{sse-z3z3}

 Consider the probability $\mu$ such that
 $\mu(a)=\mu(b)=p,\mu(a^2)=\mu(b^2)=q=1/2-p$. According to Proposition
 \ref{pr-stat}, the harmonic measure is ergodic.
 The Stationary Traffic Equations \eref{eq-trafficstat} can be solved
 explicitly. We get
\begin{equation}\label{eq-z3z3a=b}
 r(a)=r(b)=\frac{4p-3+\sqrt{16p^2-8p+5}}{4(4p-1)}, \quad
 r(a^2)=r(b^2)=\frac{4p+1-\sqrt{16p^2-8p+5}}{4(4p-1)}\:.
\end{equation}
In particular, the drift is
 \begin{equation}\label{eq-driftz3z3a=b}
 \gamma \ = \ \frac{1}{2}-p + \bigl(2p-\frac{1}{2}\bigr)r(1) \ = \ -\frac{1}{4} +
 \frac{1}{4}\sqrt{16p^2-8p+5}\:.
 \end{equation}

At last, consider the case $\mu(a)=p,\mu(a^2)=q,$ and
$\mu(b)=\mu(b^2)=(1-p-q)/2$. 
Here, it is not difficult
to check that there is no solution to the Stationary Traffic Equations.
However, solving explicitly the Traffic Equations
is feasible. The formulas for $r$ are too lengthy to be
reproduced here. But, for the drift, several simplifications occur
and we obtain the following formula: 
\[
\gamma= 2(1-p-q)\sqrt{\frac{p^2+q^2+pq}{p^2+q^2-2pq+3}}\:.
\]
In the subcase $p=q$, we get $\gamma=2p(1-2p)$, a formula that can
also be obtained using \eref{eq-elementary-general}.
For the general nearest neighbor random walk on $\Z/3\Z\star \Z/3\Z$, we
did not succeed in completely solving the Traffic Equations. 

\subsection{The simple random walks on $\Z/k\Z\star \Z/k\Z$}
\label{sse-simplezn} 

\begin{figure}[ht]
\[ \epsfxsize=200pt \epsfbox{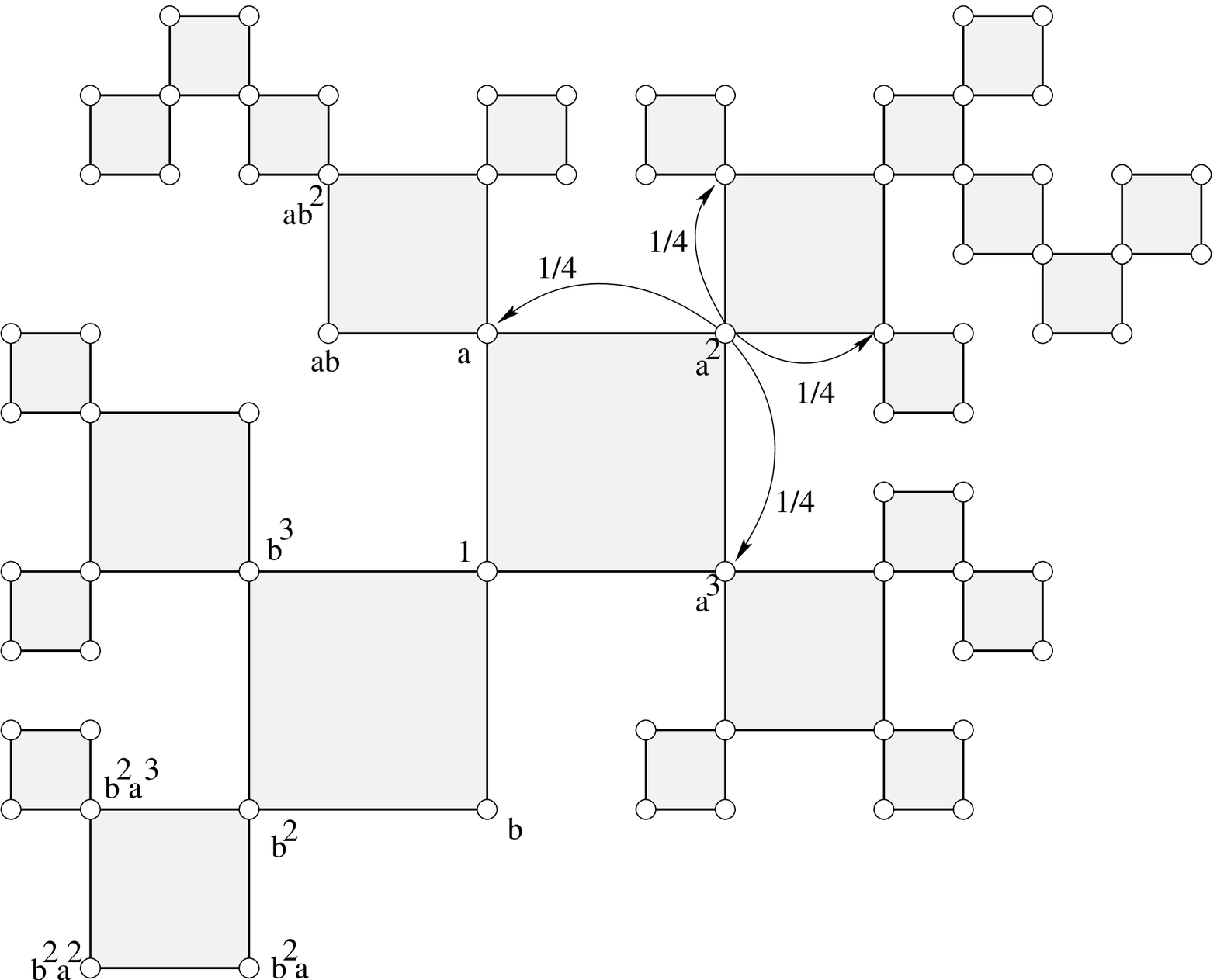} \qquad \epsfxsize=180pt \epsfbox{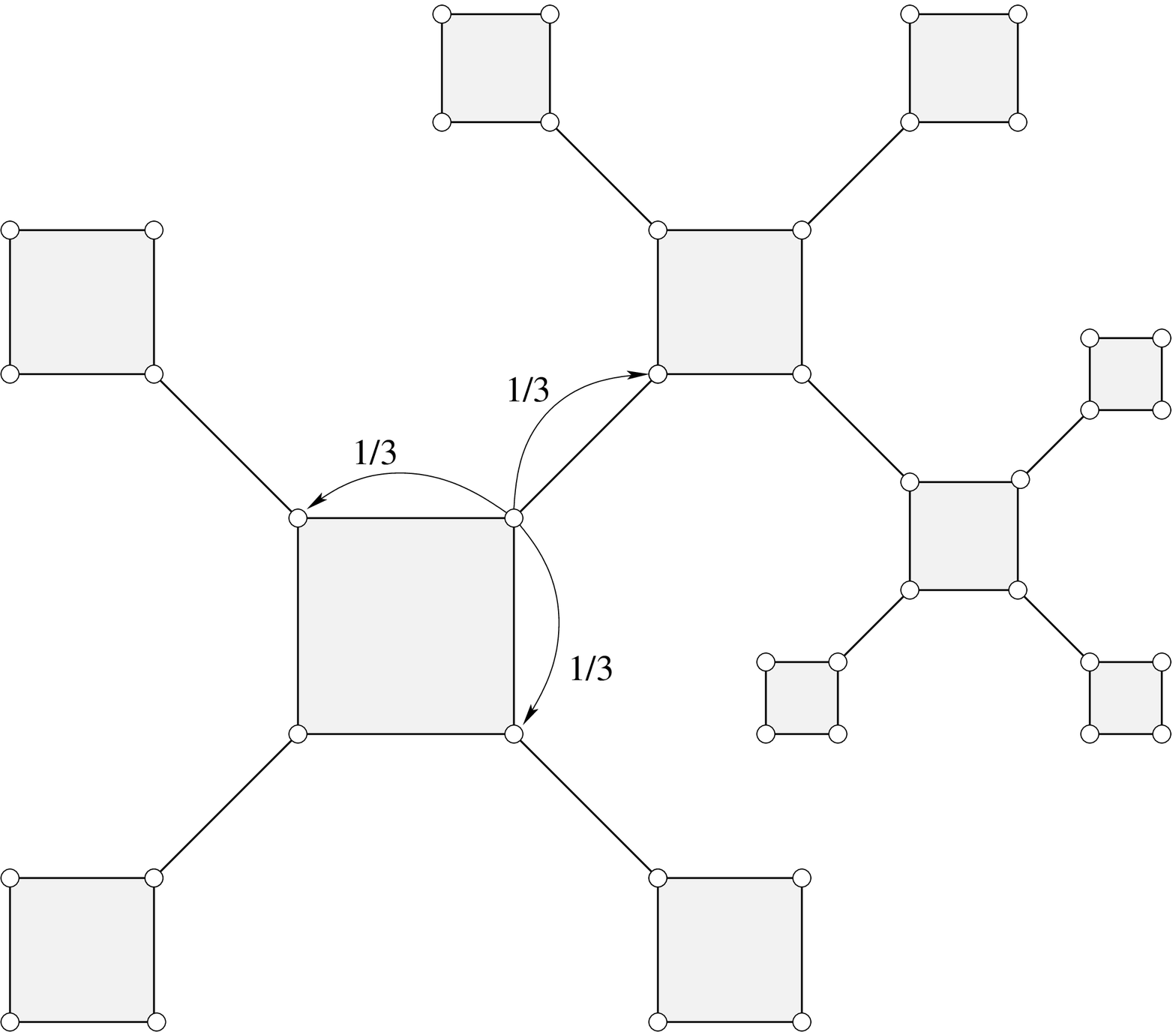}\]
\caption{The simple random walk on $\Z/4\Z\star \Z/4\Z$ (left) and
  $\Z/2\Z\star \Z/4\Z$ (right).}\label{fi-cayleyz4z4}
 \end{figure}

We now consider the whole family of groups $\Z/k\Z\star \Z/k\Z, k\geq
3$. However, we have to compromise by
considering only simple random walks with respect to a minimal set of
generators. In  Figure \ref{fi-cayleyz4z4} (left), we show this simple 
random walk in the case $\Z/4\Z\star \Z/4\Z$. 

\medskip

\medskip

We obtain ``semi-explicit'' formulas: we
define recursively a family of polynoms of one variable $(F_k)_k$, and 
the harmonic measure is expressed as a function of the 
unique solution in $(0,1)$ of $F_k=1$. For $k\geq 6$, we have no
closed form formula for this unique root. 

\medskip

Consider the free product $G_1\star G_2=\Z/k\Z\star \Z/k\Z$. 
Let $a$ and $b$ be the generators of the two cyclic groups. 
Consider the simple random walk
    $(\Z/k\Z\star \Z/k\Z,\mu)$ with
    $\mu(a)=\mu(b)=\mu(a^{-1})=\mu(b^{-1})=1/4$. 

According to Proposition \ref{pr-stat}, the harmonic measure is
ergodic. To determine it, we need to solve the Stationary
Traffic Equations \eref{eq-trafficstat}. Set
$r(i)=r(a^i)+r(b^i)=2r(a^i)=2r(b^i)$. The Stationary Traffic Equations
are: 
\begin{equation}\label{eq-master}
r= r \Bigl[ \frac{1}{4} 
\left[\begin{array}{ccccc}
0 & 1 &  & &   \\
1 & 0 & 1 &  &   \\
  & \ddots & \ddots & \ddots &   \\
&&1&0 & 1 \\
&&& 1 & 0 
\end{array}\right] + \frac{1}{4} 
\left[\begin{array}{ccccc}
1 & 0 & \cdots & 0 & 1 \\
1 & 0 & \cdots & 0 & 1 \\
\vdots & \vdots & & \vdots & \vdots \\
1 & 0 & \cdots & 0 & 1 \\
1 & 0 & \cdots & 0 & 1
\end{array}\right] + \frac{1}{4}
\left[\begin{array}{cccc}
r(1) & r(2) & \cdots & r(k-1) \\
0 & \cdots & \cdots  & 0 \\
\vdots  & &   & \vdots \\
0 & \cdots & \cdots  & 0 \\
r(1) & r(2) & \cdots & r(k-1) 
\end{array}\right] \Bigr]\:.
\end{equation}

If $[r(1),r(2),\ldots , r(k-1)]$ is a solution to \eref{eq-master}, then it is
clear that $[r(k-1),\ldots,r(2),r(1)]$ is also a solution. By
uniqueness, we deduce that  $r(i)=r(k-i)$ for all $i$. In
particular, we have $r(1)=r(k-1)$. 
For convenience, set $r(0)=r(k)=1$. Then the
Equations \eref{eq-master} rewrite as:
\begin{equation}\label{eq-sym}
\forall i \in \{1,k-1\}, \qquad r(i+1)+r(i-1)=2(2-r(1))r(i) \:.
\end{equation}

Consider the applications $F_n: [0,1]\rightarrow \R, n\in \N,$ defined by 
\begin{equation}\label{eq-rec}
F_0(x)=1, \quad F_1(x) =x, \quad \forall n\geq 2, \ F_n(x) = 2(2-x)F_{n-1}(x)
-F_{n-2}(x)\:.
\end{equation}

To illustrate, here are the first values of $F_i$:
\begin{equation*}
F_2  =  - 2 x^2 +4 x - 1, \quad 
F_3=  4x^3 -16x^2 + 17 x -4, \quad
F_4=  -8x^4 +48x^3 -96x^2 +72x-15 \:.
\end{equation*}

By construction, the unique solution $r\in \mathring{\cB}$ of the Traffic Equations
\eref{eq-sym} satisfies $r(i)=F_i(r(1))$ for all 
$i\in \{1,\ldots, k-1\}$. Hence, it is enough to determine $r(1)$. 
Set $x_k=r(1)$. We have $F_i(x_k)=F_{k-i}(x_k)$ for all $i\in
\{1,\ldots,k-1\}$. 
Now, using
\eref{eq-rec} twice, we get
\begin{eqnarray*}
F_k(x_k) & = & 2(2-x_k)F_{k-1}(x_k)
-F_{k-2}(x_k) \ = \ 2(2-x_k)F_{1}(x_k)-F_{2}(x_k) \\
& = & 2(2-x_k)F_{1}(x_k)-\bigl[
  2(2-x_k)F_{1}(x_k) - F_0(x_k)\bigr] \ = \ F_0(x_k) \ = \ 1\:.
\end{eqnarray*}

Next lemma shows that the equality $F_k(x_k)=1$ is actually a
characterization of $x_k$. 

\begin{lemma}\label{le-unique}
For $k\geq 3$, the equation $F_k(x)=1$ has a unique solution $x_k$ in
$(0,1)$. 
\end{lemma}

The results are gathered in the theorem below. 

\begin{theorem}\label{th-zkzk}
Consider the group $\Z/k\Z\star \Z/k\Z$, $k\geq 3$, the generators of the two
cyclic groups being respectively $a$ and $b$. Consider the simple
random walk $(\Z/k\Z\star \Z/k\Z,\mu)$ with
$\mu(a)=\mu(a^{-1})=\mu(b)=\mu(b^{-1})=1/4$. 
The harmonic measure is
the ergodic Markovian multiplicative measure associated with 
\[
r= [ x_k, F_2(x_k), F_3(x_k),\ldots , F_2(x_k), x_k ]
\]
where $x_k$ is the unique solution in $(0,1)$ of the equation
$F_k(x)=1$. The drift is $\gamma_k= (1-x_k)/2$. It is a strictly increasing function of $k$ and
$\lim_k \gamma_k=1/3$. 
\end{theorem}

To illustrate,
here is the vector $r$ for $\Z/4\Z\star \Z/4\Z$:
\begin{equation}\label{eq-rz4z4}
\bigl[r(a),r(a^2),r(a^3)\bigr]=\bigl[r(b),r(b^2),r(b^3)\bigr]= \Bigl[
  \frac{3-\sqrt{5}}{8}, \frac{\sqrt{5}}{4} -\frac{1}{2}, \frac{3-\sqrt{5}}{8} \Bigr]\:.
\end{equation}
Now, here is a table of the first values of $\gamma_k$, given
either in closed form or numerically when no closed form could be
found. Set $\Z_k=\Z/k\Z$. 

\begin{center}
\begin{tabular}{|l|c|c|c|c|c|c|} \hline
 & $\Z_3\star \Z_3$ & $\Z_4\star \Z_4$ &   $\Z_5\star
 \Z_5$ &  $\Z_6\star \Z_6$ & $\Z_7\star \Z_7$ & $\Z_8\star
 \Z_8$ \\ \hline 
&&&&&& \\
$\gamma$ &  $1/4$ & $(\sqrt{5}-1)/4$  & $(\sqrt{13}-1)/8$  &
  0.330851... & 0.332515...  &  0.333062...\\ 
&&&&&& \\\hline
\end{tabular}
\end{center}

\subsection{The simple random walks on $\Z/2\Z \star \Z/k\Z$}
\label{sse-z2zk}

We study simple random walks on
$\Z/2\Z \star \Z/k\Z$ for a minimal set of generators. 
The model considered is illustrated in Figure \ref{fi-cayleyz4z4}
(right) in the case $\Z/2\Z \star \Z/4\Z$. 
The approach and results are similar to the ones in \S \ref{sse-simplezn}.
The groups
$\Z/2\Z\star \Z/k\Z$ are known as the Hecke groups. 

\medskip


Consider the applications $G_n: [0,1]\rightarrow \R, n\in \N,$ defined by 
\begin{equation}\label{eq-recz2zk}
G_0(x)=\frac{1}{4}+\frac{x}{2}, \quad G_1(x) =x, \quad \forall n\geq
2, \ G_n(x) = \frac{8(1-x)}{3-2x}G_{n-1}(x) 
-G_{n-2}(x)\:.
\end{equation}

\begin{theorem}\label{th-z2zk}
Consider the group $\Z/2\Z\star \Z/k\Z$, $k\geq 3$, the generators of the two
cyclic groups being respectively $a$ and $b$. Consider the simple
random walk $(\Z/2\Z\star \Z/k\Z,\mu)$ with
$\mu(a)=\mu(b)=\mu(b^{-1})=1/3$. 
Let the functions $G_i$ be defined as in \eref{eq-recz2zk} and let
$y_k$ be the
unique solution in $(0,1/2)$ of $G_{k-1}(y_k)=y_k$. 
The harmonic measure is
the Markovian multiplicative measure associated with $r$: $r(a) =
G_0(y_k), \ \forall i
\in \{1,\ldots, k-1\}, \ r(b^i)=G_i(y_k)$. The drift is $\gamma_k=
(1-2y_k)/3$. It is a strictly increasing function of $k$ and $\lim_k
\gamma_k =2/9$. 
\end{theorem}

Here is the vector $r$ for $\Z/2\Z\star \Z/4\Z$:
\[
\bigl[r(a),r(b),r(b^2),r(b^3)\bigr]= \Bigl[
\frac{7-\sqrt{7}}{12}, \frac{2}{3} - \frac{\sqrt{7}}{6},
\frac{-11+5\sqrt{7}}{12} , \frac{2}{3} - \frac{\sqrt{7}}{6}\Bigr]\:. 
\]
Here is the drift, either in closed form or numerically, for the small values
of $k$. 

\begin{center}
\begin{tabular}{|l|c|c|c|c|c|c|} \hline
& $\Z_2\star \Z_3$ & $\Z_2\star \Z_4$ &   $\Z_2\star
\Z_5$ &  $\Z_2\star \Z_6$ & $\Z_2\star \Z_7$ & $\Z_2\star
\Z_8$ \\ \hline 
 &&&&&& \\
 $\gamma$ &  $2/15$ & $(\sqrt{7}-1)/9$  & $(2\sqrt{61}-4)/57$  &
   0.213412... & 0.217921...  &  0.220101...\\ 
 &&&&&& \\\hline
 \end{tabular}
 \end{center}

\section{Entropy and Extremal Generators}\label{se-vershik}

We study Vershik notion of extremal generators for
free products of finite groups. The results are of two kinds,
structural results derived from the shape of the harmonic measure
(Prop. \ref{pr-extr} and \ref{Hausdorff}), and
explicit computations obtained by solving the Traffic Equations
(Section \ref{sse-explicitentropy}). 

\medskip

The {\em entropy} of a probability measure $\mu$ with finite support $S$ is
defined by
\begin{equation}\label{eq-entropy}
H(\mu) = -\sum_{x\in S} \mu(x) \log [\mu(x) ] \:.
\end{equation}
The entropy is maximized for the uniform measure $U$ on $S$.

Let $G$ be a free product of finite groups and let 
$(G,\mu)$ be a random walk defined as in \S \ref{se-prel}. 
Let $(X_n)_n$ be a realization of the random walk. 
The {\em entropy of $(G,\mu)$}, introduced by Avez \cite{avez}, is 
\begin{equation}\label{eq-avez}
h =  h(\mu)=\lim_n \frac{H(\mu^{*n})}{n}= \lim_n -\frac{1}{n}
\log\mu^{*n}(X_n) \:,
\end{equation}
$a.s.$ and in $L^p$, for all $1\leq p <\infty$. The existence of the
limits as well as their equality follow from Kingman's Subadditive
Ergodic Theorem (see \cite{avez} and Derriennic \cite{derr}). In
general, the entropy is difficult to compute, see \cite{derr}. But 
not in our case. 

\begin{proposition}
Let $r$ be the unique solution in $\mathring{\cB}$ of the Traffic
Equations. Set $q(a)=r(a)/r(\Sigma\moins \Sigma_a)$ for all $a\in
\Sigma$. The entropy is given by:
\begin{eqnarray}
h & = & - \sum_{a\in \Sigma} \mu(a) \int \log \bigl[
  \frac{da^{-1}\mu^{\infty}}{d\mu^{\infty}}(\xi)\bigr]
  d\mu^{\infty}(\xi) \label{eq-entro1}\\
  & = & -  \sum_{a\in \Sigma} \mu(a) \Bigl[ \log \bigl[
  \frac{1}{q(a^{-1})} \bigr]
  r(a^{-1}) + \sum_{b\in \Sigma_{a}\moins a^{-1}} \log \bigl[
  \frac{q(ab)}{q(b)}\bigr] r(b) + \log[q(a)]\sum_{b\in \Sigma \moins
  \Sigma_a} r(b) \Bigr]  \label{eq-entro2}
\end{eqnarray}
where $da^{-1}\mu^{\infty}/d\mu^{\infty}$ is the Radon-Nikodym derivative of $a^{-1}\mu^{\infty}$ with
respect to $\mu^{\infty}$. 
\end{proposition}

Formula \eref{eq-entro1} is proved for instance 
in \cite[Corollary 4.5]{ledr00} in the context of the free group. 
The proof adapts easily to the present setting. 
The formulation in \eref{eq-entro2} follows directly from the one in
\eref{eq-entro1}, using that the harmonic measure is Markovian
multiplicative. 

\medskip

The {\em volume} of the group $G$ with respect to the finite 
set of generators $\Sigma$ is 
\begin{equation}\label{eq-volume}
v = v(\Sigma) = \lim_n \frac{1}{n} \log \#\{g\in G, |g|=n\}\:.
\end{equation}
The limit exists by subadditivity. The following
fundamental inequality was proved and
highlighted in
\cite{vers00} (see also \cite{guiv80} and \cite{kaim98,ledr83}): 
\begin{equation}\label{eq-fundineq}
h/\gamma \leq v\:,
\end{equation}
where we recall that $\gamma$ is the drift of the random walk. The
interpretation is that the proportion of typical elements visited by
the walk is less than or equal to the total number of elements. 

\medskip

This inequality is reminiscent of the one between metric entropy and
topological entropy in symbolic dynamics \cite{LiMa}. 
Also, for a free product of finite groups, $h/\gamma$ and $v$ can be 
interpreted respectively as the Hausdorff dimension of the harmonic measure $\mu^{\infty}$ 
and the Hausdorff dimension of its support. See Prop. \ref{Hausdorff} for a precise statement.  

\medskip

Observe
that $v$ and $\gamma$ depend on $\Sigma$ but not $h$, and that $h$ and
$\gamma$ depend on $\mu$ but not $v$. 
It might be enlightning to view the inequality as
$h(\mu)/v(\Sigma)\leq \gamma(\Sigma,\mu)$. 
Define 
\[
Q(\Sigma)= \sup_{\mu\in {\cal S}} \frac{
h(\mu)}{\gamma(\Sigma,\mu)v(\Sigma) }\:,
\]
where $\cS$ 
is the set of probability measures on $\Sigma$ such that: (i)
$\forall a, \mu(a)=\mu(a^{-1})$ ($\mu$ is symmetric), and (ii) $\bigcup_n
\mbox{supp} (\mu^{*n})= G$. The set $\cS$ is not closed due to
condition (ii). The closure $\overline{\cS}$ is the set of probability
measures on $\Sigma$ satisfying only (i). 
It might be the case that the supremum in $Q(\Sigma)$ is attained only for a measure
$\mu\in \overline{\cS}\moins \cS$, for which the random walk is not
transient, and for which $\gamma=0$.  
We will see an occurence of this situation in Section \ref{sse-explicitentropy}. 

\medskip

It is suggested by Vershik \cite{vers00} to view $Q(\Sigma)$
as a measure of the `quality' of
the set of generators $\Sigma$. If $Q(\Sigma)=1$ then $\Sigma$ is said to be
{\em extremal}. 

\medskip

This notion is investigated in \cite{VNBi} for braid
groups and  partially commutative free groups. In particular,
it is proved that the natural generators
are not extremal for partially commutative free groups of the type
$\pres{a_1,\ldots , a_n}{\forall i=1,\ldots, n-1, \
  a_ia_{i+1}=a_{i+1}a_i}$, for $n$ large enough. 
Below, we continue the study of this notion using free products of
finite groups. 

\medskip

Consider a free product $G_1\star
\cdots \star G_n$ with $|G_i|=|G_j|$ for all $i,j$.  
Set $\Sigma =\sqcup_i G_i\moins\{1_{G_i}\}$. It is elementary to prove that $\Sigma$ is
extremal. Indeed consider the
uniform distribution $\mu$ on $\Sigma$. It is clear by symmetry that 
the harmonic measure should be the
uniform measure on $L^{\infty}$, that is, $\forall u,v \in L\cap\Sigma^k$,
$\mu^{\infty}(u\Sigma^{\N})= \mu^{\infty}(v\Sigma^{\N})$. The exact
values of $h,\gamma,$ and $v$ follow readily. Set $k=|G_i\moins \{1_{G_i}\}|$, and
$K=nk=|\Sigma|$. We have
\[
v=\log\bigl[K-k\bigr], \quad \gamma= -
\frac{1}{K} + (1-\frac{k}{K}), \quad h
=\log\bigl[K-k \bigr] \bigl(-
\frac{1}{K} + (1-\frac{k}{K})\bigr) \:.
\]
In particular, $h/\gamma=v$.  This is mentioned in \cite[Section 2.5]{vers00}.

The argument collapses when the different groups $G_i$ do not have the
same cardinality. In particular, there is no reason to expect the
uniform ditribution on $\Sigma$ to maximize $h/(\gamma v)$.
Indeed, the result $Q(\Sigma)=1$ is still true but non-elementary.
We adress this point in Proposition \ref{pr-extr}. 
For definitions and details regarding subshifts and measures
of maximal entropy, see \cite{LiMa}.

\begin{proposition}\label{pr-extr}
Consider a free product of finite groups $G_1\star
\cdots \star G_n$. Set $\Sigma_i=G_i\moins \{1_{G_i}\}$ and set
$k_i=|\Sigma_i|$. Let $\rho$ be the unique positive
solution of the equation: $\sum_{i=1}^n k_i/(x+k_i) =1$. Consider the
probability measure $\mu$ on $\Sigma=\sqcup_i \Sigma_i$ defined by
$\mu =\sum_i k_i/(\rho+k_i)U_{\Sigma_i}$, where $U_{\Sigma_i}$ is the
uniform distribution on $\Sigma_i$. The following properties hold. 
\par
\quad (i) We have $h/(\gamma v)=1$ for $(G,\mu)$. The generators
$\Sigma$ are extremal. 
\par
\quad (ii) The harmonic measure $\mu^{\infty}$ associated with $(G,\mu)$
is defined by

\begin{equation}\label{eq-harmonic}
\forall u_1\cdots u_k \in L, u_k \in \Sigma_j, \quad
\mu^{\infty}(u_1\cdots u_k\Sigma^{\N})= \frac{1}{\rho^{k-1}}
\frac{1}{(\rho+k_j)}\:.
\end{equation}

\quad (iii) Consider the shift: $\tau:\Sigma^{\N}\rightarrow
\Sigma^{\N}, \ (x_n)_n \mapsto (x_{n+1})_n$. Then the symbolic dynamical system $(L^{\infty},\tau)$ 
is a subshift of finite type whose  measure of maximal entropy $\nu_{max}$ is defined by

\begin{equation}\label{eq-max-entropy}
\forall u_1\cdots u_k \in L, u_1 \in \Sigma_i, u_k \in \Sigma_j,\quad
\nu_{max}(u_1\cdots u_k\Sigma^{\N})= \frac{1}{(\rho+k_i)} \frac{1}{\rho^{k-2}} \frac{1}{(\rho+k_j)}\:.
\end{equation}

\end{proposition}

\begin{proof}
Let us compute the volume $v$ of $G=(G_1\star \cdots \star
G_n)$. Consider the matrix
\begin{equation}\label{eq-A}
A= \left[ \begin{array}{ccccc}
0 & k_2 & k_3 & \cdots & k_n \\
k_1 & 0 & k_3 & \cdots & k_n \\
k_1 & k_2 & 0 & \cdots & k_n \\
\vdots &\vdots &&&\vdots \\
k_1 & k_2 & \cdots & 0 & k_n \\
k_1 & k_2 & \cdots  & k_{n-1} & 0
\end{array}\right]\:.
\end{equation}
Clearly $\sum_{ij} k_iA^{n-1}_{ij}$ is equal to the number of elements
of $G$ of length $n$ (with respect to $\Sigma$). It follows that
$v=\log(\rho)$, where
$\rho$ is the spectral radius of $A$. According to the
Perron-Frobenius Theorem, there exists a unique $Y=[Y_1,\ldots, Y_n], Y_i >0,
\sum_i Y_i=1$, such that $YA=\rho Y$, i.e. $\forall i$,
\begin{equation*}
k_i (\sum_{j\neq i} Y_j) = \rho Y_i \ \iff \ k_i(1-Y_i) =\rho Y_i \
\iff \ Y_i
= \frac{k_i}{\rho + k_i}\:.
\end{equation*}
Therefore, the spectral radius $\rho$ satisfies:
\begin{equation}\label{eq-rho}
\frac{k_1}{\rho + k_1} + \cdots + \frac{k_n}{\rho + k_n} =1\:.
\end{equation}
Since the homographic maps $k_i/(x+k_i)$ are strictly decreasing
functions of $x$ in $\R_+$, the Equation \eref{eq-rho} characterizes
$\rho$.

Now consider the
probability measure $\mu$ on $\Sigma$ defined by
$\mu =\sum_i k_i/(\rho + k_i)U_{\Sigma_i}$, where $U_{\Sigma_i}$ is the
uniform distribution on $\Sigma_i$. For simplicity, set
$\mu_i = \mu(\Sigma_i)/k_i= 1/(\rho + k_i)=\mu(a)$  for all $a\in
\Sigma_i$. Let $r$ be the
solution to the Traffic Equations and let $q$ be the probabilities of
ever hitting the generators.
Set $r_i = r(\Sigma_i)/k_i=
r(a), \forall a\in \Sigma_i$, and $q_i = q(\Sigma_i)/k_i=
q(a), \forall a\in \Sigma_i$. Observe that we have
\begin{equation}\label{eq-interm}
\sum_i k_i\mu_i =1 , \qquad \sum_i k_ir_i =1 , \qquad \forall i, \ q_i =
\frac{r_i}{\sum_{j\neq i} k_jr_j} = \frac{r_i}{1-k_ir_i} \:.
\end{equation}
Using Formulas \eref{eq-drift2} and \eref{eq-entro2}, we get the drift:
\begin{equation}\label{eq-dri}
\gamma = \sum_i k_i r_i \ \bigl[ -\mu_i +
  \sum_{j\neq i} k_j\mu_j \bigr] = 1 -\sum_i k_i(k_i+1)\mu_i r_i  \:,
\end{equation}
and the entropy:
\begin{equation}\label{eq-ent}
h= -\sum_i \mu_i \log(q_i) \bigl[ -r_i +\sum_{j\neq i} k_jr_j
  \bigr] = \sum_i \mu_i \log(1/q_i) - \sum_i \mu_i \log(1/q_i)r_i
  \bigl[ k_i +1 \bigr] \:.
\end{equation}
To obtain the formula for $h$, we have used that for all $i$ and all
$a,b\in \Sigma_i$, such that $a\ast b\in \Sigma_i$, we have
$q(a)=q(a\ast b)$, so $\log[q(a\ast b)/q(a)]=0$.

Assume that we have: $\forall i, \log(1/q_i)=v$. Then it follows from
\eref{eq-dri} and \eref{eq-ent} that $h/\gamma=v$.
Let us prove that $\log(1/q_i)=v$ for all $i$.
Recall that $q_i$ satisfies the Equations \eref{eq-reach}:
\begin{equation}\label{eq-qi}
q_i = \mu_i + (k_i\mu_i - \mu_i)q_i + q_i \sum_{j\neq i} k_j\mu_j
q_j \:.
\end{equation}
Assume that $q_i=q_j=q$ for all $i,j$. Then we get
\[
(1-k_i\mu_i) q^2 - (1-k_i\mu_i + \mu_i)q +
\mu_i = 0 \:.
\]
The two solutions of this second order equation are:
$q =1$ and $q = \mu_i/(1-k_i\mu_i)$.
Since we must have $0<q<1$, the right solution is the second one, and we deduce that:
\[
\mu_i = \frac{1}{1/q + k_i } \quad \implies \quad \sum_i
\frac{k_i}{1/q +k_i} = 1 \:.
\]
Comparing with \eref{eq-rho}, we conclude that
$1/q=\rho$. Backtracking in the argument, we conclude that $\forall i,
q_i = 1/\rho,$ is  a solution to the Equations \eref{eq-qi}.
Now the Equations \eref{eq-qi} have at most one solution in
$(0,1)^{\Sigma}$, otherwise, following the proof of Theorem
\ref{th-jackpot}, one would get several solutions in $\mathring{\cB}$ for
the Traffic Equations, which is impossible.
We conclude that: $\forall i, \log(1/q_i)=\log(\rho)= v$. It
completes the proof of $h/(\gamma v)=1$.

\medskip

Observe that: $r_i = \mu_i = Y_i/k_i = 1/(\rho + k_i)$.
Applying Theorem \ref{th-jackpot}, we get that the harmonic measure
$\mu^{\infty}$ is given by \eref{eq-harmonic}.

\medskip

Consider the translation shift: $\tau:\Sigma^{\N}\rightarrow
\Sigma^{\N}, \ (x_n)_n \mapsto (x_{n+1})_n$.
The symbolic dynamical system $(L^{\infty},\tau)$ is clearly 
a subshift of finite type. For such a system, the measure of
{\em maximal entropy} is well-known, and it turns out to be precisely
the measure $\nu_{\max}$ given in \eref{eq-max-entropy}.
\end{proof}

We now turn our attention to Hausdorff dimensions. The Hausdorff
dimensions of a metric space $(X,d)$ and a Borel measure $\nu$
on $X$ are denoted respectively by $\text{HD}(X)$ and $\text{HD}(\nu)$,
see for instance \cite{matt95} for the definitions.

\begin{proposition}\label{Hausdorff}
Let $G=G_1\star \cdots \star G_n$ be a free product of finite groups. Set
$\Sigma_i= G_i\moins\{1_{G_i}\}$ and $\Sigma = \sqcup_i \Sigma_i$. 
Let $v$ be the volume of $G$ with respect to $\Sigma$. 
Denote by $d$ the metric defined on the set $L^{\infty}$ of right-infinite
normal form words by $d (\xi_{1},\xi_{2})\,=\,e^{-|\xi_{1} \wedge \xi_{2}|}$ 
where $\xi_{1} \wedge \xi_{2}$ is the greatest common prefix of
$\xi_{1}$ and $\xi_{2}$.  Let $\nu_{max}$ be the measure of maximal entropy of the
dynamical system $(L^{\infty},\tau)$, see \eref{eq-max-entropy}. The
following properties hold. 
\par
\quad (i) We have $\text{{\rm HD}}(\nu_{max})\:=\:\text{{\rm
    HD}}(L^{\infty}) \:=\: v$. 
\par
\quad (ii) Let $\mu$ be a probability on $\Sigma$ whose support generates $G$. 
We have $\text{{\rm HD}}(\mu^{\infty})\:=\: h/\gamma$.
\end{proposition}

\begin{proof}
The equality $\text{{\rm
    HD}}(L^{\infty}) \:=\: v$ is analogous to
\cite[Proposition 1.9]{ledr00}. 
By taking some precautions, the proof can be adapted. Let $\mu^{\infty}$ be the harmonic measure
associated with the probability $\mu$ defined in Prop.
\ref{pr-extr}. With $\mu^{\infty}$ in the role of the measure
$m$, the proof of \cite[Proposition 1.9]{ledr00}
goes through. Here, it is central to know that $\mu^{\infty}$ is
    closely related to $\nu_{\max}$ to be able to adapt
the argument from \cite[Proposition 1.9]{ledr00}. The equality $\text{{\rm
    HD}}(\mu^{\infty})\:=\: h/\gamma$ is proved in \cite[Theorem 4.15]{ledr00} for a finitely generated free
group. The proof adapts easily to the present setting.
\end{proof}

\subsection{Explicit computations}\label{sse-explicitentropy} 

Let us go back to extremal generators. After having settled Prop.
\ref{pr-extr}, other natural questions become in order.
Is there a unique measure on $\Sigma$ maximizing $h/(\gamma
v)$? Are other sets of generators extremal?

\medskip

We believe that some new light is shed on these questions by the three
results below. 
These results are obtained by solving explicitly
the Traffic Equations in the cases considered. 
The proofs and the details are to be found in the Math ArXiv
appendix \cite{MaMa05}.

\medskip

\hspace*{1cm} $\bullet$ A. Consider a free product of two finite
groups, $G_1\star G_2$. 
Consider the probability $\mu_p$ on $\Sigma$ such
that $\mu_p = p U_{\Sigma_1} + (1-p) U_{\Sigma_2}$,
where $U_{\Sigma_i}$ is the uniform distribution on $\Sigma_i$, and
where $p\in (0,1)$. Then we have
$h/(\gamma v)=1$, for the whole family of probabilities $\mu_p$. This
does not extend to the free product of more than 2 groups, see B,
nor to minimal sets of generators $S \subset \Sigma$, see C.
According to Prop. \ref{Hausdorff}, we have $\text{{\rm
    HD}}(\mu_{p}^{\infty})\:=\:\text{{\rm HD}}(L^{\infty})$ 
for $p\in (0,1)$. This contrasts sharply with \cite[Theorem 2.1]{ledr00}.

\par

\hspace*{1cm} $\bullet$ B. Consider the group $\Z/2\Z \star \Z/2\Z \star \Z/2\Z$,
and let $a,b,$ and $c$ be the non-identity elements of the three
cyclic groups. Consider the family of probability measures $\mu_p,
p\in (0,1/2),$ defined by: $\mu_p(a)=\mu_p(b)=p, \mu_p(c)=1-2p$. Among
them, the
only probability such that $h/(\gamma v)=1$ is $\mu_{1/3}$. 

\par

\hspace*{1cm} $\bullet$ C. We now give three examples of various behaviours
for a generating set $S$ of a free product $G_1\star G_2$ such that
$S\subsetneq \Sigma$, the natural generators. 
\par
\qquad First, consider the group $\Z/2\Z \star \Z/4\Z$. Let $a$ and
$b$ be respective generators of the two cyclic groups. The
minimal set of generators $S=\{a,b,b^{-1}\}$ is extremal but $h/(\gamma v)<1$ 
for all $\mu$ in $\cS$. 
\par
\qquad Consider now the group $\Z/3\Z \star \Z/4\Z$. The minimal set
of generators $S=\{a,a^{-1},b,b^{-1}\}$ is extremal. Actually, the only
symmetric probability measure on $S$ for which $h=\gamma v$ is
$\mu \,=\, p\,(\delta_{a}+\delta_{a^-1})+(1/2-p)(\delta_{b}+\delta_{b^-1})$ 
where $p\,=\,0.432693\cdots$ is the middle root of the
polynomial $5x^3-13x^2+7x-1$. Observe that $\mu \in \cS$.
\par
\qquad Last, consider the group $\Z/4\Z \star \Z/4\Z$. The minimal set
of generators $S=\{a,a^{-1},b,b^{-1}\}$ is not extremal. Indeed, 
\[
Q(S)=  \frac{5+\sqrt{5}}{4}\frac{\log \bigl(1/2 + \sqrt{5}/2\bigr)}
{\log(1+\sqrt{2})} = 0.987686\cdots \:.
\]

\paragraph{Acknowledgement.} The authors would like to thank 
the following colleagues for stimulating discussions: Christian
Blanchet, Philippe Bougerol, Thierry Bousch, Xavier Bressaud, Yves Guivarc'h,
Fran\c{c}ois Ledrappier, Vincent Le Prince, James Martin, Sergei Nechaev, Vadim
Kaimanovich, and Wolfgang Woess.


\begin{thebibliography}{10}

\bibitem{avez}
A.~Avez.
\newblock Entropie des groupes de type fini.
\newblock {\em C. R. Acad. Sci. Paris S\'er. A-B}, 275:1363--1366, 1972.

\bibitem{cart88}
D.~Cartwright.
\newblock Some examples of random walks on free products of discrete groups.
\newblock {\em Ann. Mat. Pura Appl. (4)}, 151:1--15, 1988.

\bibitem{CaSo86}
D.~Cartwright and P.~Soardi.
\newblock Random walks on free products, quotients and amalgams.
\newblock {\em Nagoya Math. J.}, 102:163--180, 1986.

\bibitem{derr}
Y.~Derriennic.
\newblock Quelques applications du th\'eor\`eme ergodique sous-additif.
\newblock {\em Ast\'erisque}, 74:183--201, 1980.

\bibitem{DyMa}
E.~Dynkin and M.~Malyutov.
\newblock Random walk on groups with a finite number of generators.
\newblock {\em Sov. Math. Dokl.}, 2:399--402, 1961.

\bibitem{furs}
H.~Furstenberg.
\newblock Noncommuting random products.
\newblock {\em Trans. Amer. Math. Soc.}, 108:377--428, 1963.

\bibitem{furs71}
H.~Furstenberg.
\newblock Random walks and discrete subgroups of {L}ie groups.
\newblock {\em Advances Probab. Related Topics}, 1:1--63, 1971.

\bibitem{GeWo}
P.~Gerl and W.~Woess.
\newblock Local limits and harmonic functions for nonisotropic random walks on
  free groups.
\newblock {\em Probab. Theory Relat. Fields}, 71(3):341--355, 1986.

\bibitem{guiv80}
Y.~Guivarc'h.
\newblock Sur la loi des grands nombres et le rayon spectral d'une marche
  al\'eatoire.
\newblock {\em Ast\'erisque}, 74:47--98, 1980.

\bibitem{kaim98}
V.~Kaimanovich.
\newblock Hausdorff dimension of the harmonic measure on trees.
\newblock {\em Ergodic Theory Dynam. Systems}, 18(3):631--660, 1998.

\bibitem{kaim00}
V.~Kaimanovich.
\newblock The {P}oisson formula for groups with hyperbolic properties.
\newblock {\em Ann. of Math. (2)}, 152(3):659--692, 2000.

\bibitem{king}
J.~Kingman.
\newblock Subadditive ergodic theory.
\newblock {\em Annals of Probability}, 1:883--909, 1973.

\bibitem{lall93}
S.~Lalley.
\newblock Finite range random walk on free groups and homogeneous trees.
\newblock {\em Ann. Probab.}, 21(4):2087--2130, 1993.

\bibitem{ledr83}
F.~Ledrappier.
\newblock Une relation entre entropie, dimension et exposant pour certaines
  marches al\'eatoires.
\newblock {\em C. R. Acad. Sci. Paris S\'er. I.}, 296:369--372, 1983.

\bibitem{ledr00}
F.~Ledrappier.
\newblock Some asymptotic properties of random walks on free groups.
\newblock In J.~Taylor, editor, {\em Topics in probability and Lie groups:
  boundary theory}, number~28 in CRM Proc. Lect. Notes, pages 117--152.
  American Mathematical Society, 2001.

\bibitem{LiMa}
D.~Lind and B.~Marcus.
\newblock {\em An introduction to symbolic dynamics and coding}.
\newblock Cambridge University Press, Cambridge, 1995.

\bibitem{mair04}
J.~Mairesse.
\newblock Random walks on groups and monoids with a {M}arkovian harmonic
  measure.
\newblock {LIAFA} research report 2004-05, Universit\'e Paris 7, 2004.

\bibitem{MaMa04}
J.~Mairesse and F.~Math\'eus.
\newblock Random walks on groups with a tree-like {C}ayley graph.
\newblock In 
  {\em Mathematics and computer science. III. Algorithms, trees, combinatorics
  and probabilities}, Trends in Mathematics, pages 445--460. Birkhauser Verlag,
  2004.

\bibitem{MaMa05}
J.~Mairesse and F.~Math\'eus.
\newblock Appendix to the paper ``Random walks on free products of cyclic groups''.
\newblock arXiv:math.PR/0509208, 2005.

\bibitem{matt95}
P.~Mattila.
\newblock {\em Geometry of sets and measures in {E}uclidean spaces}, volume~44
  of {\em Cambridge Studies in Advanced Mathematics}.
\newblock Cambridge University Press, Cambridge, 1995.

\bibitem{NaWo}
T.~Nagnibeda and W.~Woess.
\newblock Random walks on trees with finitely many cone types.
\newblock {\em J. Theoret. Probab.}, 15(2):383--422, 2002.

\bibitem{nech}
S.~Nechaev.
\newblock {\em Statistics of knots and entangled random walks}.
\newblock World Scientific Publishing Co. Inc., River Edge, NJ, 1996.

\bibitem{NeVo}
S.~Nechaev and R.~Voituriez.
\newblock Random walks on three-strand braids and on related hyperbolic groups.
\newblock {\em J. Phys. A}, 36(1):43--66, 2003.


\bibitem{SaSt}
S.~Sawyer and T.~Steger.
\newblock The rate of escape for anisotropic random walks in a tree.
\newblock {\em Probab. Theory Related Fields}, 76(2):207--230, 1987.

\bibitem{taka97}
C.~Takacs.
\newblock Random walk on periodic trees.
\newblock {\em Electron. J. Probab.}, 2:no.\ 1, 1--16, 1997.

\bibitem{vers00}
A.~Vershik.
\newblock Dynamic theory of growth in groups: Entropy, boundaries, examples.
\newblock {\em Russ. Math. Surv.}, 55(4):667--733, 2000.
\newblock Translation from Usp. Mat. Nauk 55(4):59-128, 2000.

\bibitem{VNBi}
A.~Vershik, S.~Nechaev, and R.~Bikbov.
\newblock Statistical properties of locally free groups with applications to
  braid groups and growth of random heaps.
\newblock {\em Commun. Math. Phys.}, 212(2):469--501, 2000.

\bibitem{woes86b}
W.~Woess.
\newblock A description of the {M}artin boundary for nearest neighbour random
  walks on free products.
\newblock In {\em Probability measures on groups VIII}, volume 1210 of {\em
  Lecture Notes in Math.}, pages 203--215. Springer, Berlin, 1986.

\bibitem{woes86}
W.~Woess.
\newblock Nearest neighbour random walks on free products of discrete groups.
\newblock {\em Boll. Un. Mat. Ital. B (6)}, 5(3):961--982, 1986.

\bibitem{woes}
W.~Woess.
\newblock {\em Random walks on infinite graphs and groups}.
\newblock Number 138 in Cambridge Tracts in Mathematics. Cambridge University
  Press, 2000.

\end{thebibliography}

\end{document}